\documentclass[journal]{IEEEtran}

\newcommand{\beq}{\begin{equation}}
\newcommand{\eeq}{\end{equation}}
\newcommand{\bsq}{\begin{subequations}}
	\newcommand{\esq}{\end{subequations}}
\newcommand{\bq}{\begin{eqnarray}}
\newcommand{\eq}{\end{eqnarray}}
\newcommand{\bqn}{\begin{eqnarray*}}
	\newcommand{\eqn}{\end{eqnarray*}}
\usepackage{bbm}
\usepackage[pdftex]{graphicx}
\graphicspath{{Figs/}}
\usepackage{enumerate}
\usepackage{enumitem} 

\usepackage{mathptmx}       % selects Times Roman as basic font
\DeclareMathAlphabet{\mathcal}{OMS}{cmsy}{m}{n}
\usepackage{helvet}         % selects Helvetica as sans-serif font
\usepackage{courier}        % selects Courier as typewriter font
\usepackage{type1cm}        % activate if the above 3 fonts are
\usepackage{makeidx}         % allows index generation
\usepackage{graphicx}        % standard LaTeX graphics tool

\usepackage{multicol}        % used for the two-column index
\usepackage[bottom]{footmisc}% places footnotes at page bottom
\usepackage{ifthen}
\usepackage{subfigure}
\usepackage[usenames,dvipsnames]{color}

\makeindex             % used for the subject index
\usepackage{graphics} % for pdf, bitmapped graphics files
\usepackage{graphicx} % for pdf, bitmapped graphics files
\usepackage{epsfig} % for postscript graphics files
\usepackage{epstopdf}
\usepackage{mathptmx} % assumes new font selection scheme installed
\usepackage{times} % assumes new font selection scheme installed
\usepackage[cmex10]{amsmath}

\usepackage{mathtools}  % loads amsmath
\usepackage{amssymb}  % assumes amsmath package installed
\usepackage{amsthm}

\usepackage{amsfonts}
\usepackage{cite}
\usepackage{verbatim}
\usepackage{comment}
\usepackage{algorithmic}
\usepackage{multirow}
\usepackage{cite}
\usepackage{stfloats}
\usepackage{supertabular}
\usepackage{longtable}

\DeclareGraphicsExtensions{.pdf,.png,.jpg,.eps,.ps}
\renewcommand{\arraystretch}{1.2}
\usepackage{arydshln}
\usepackage{multicol}
\usepackage{colortbl}
\theoremstyle{definition}
\newtheorem{theorem}{Theorem}
\newtheorem{lemma}{Lemma}

\newtheorem{proposition}{Proposition}

\theoremstyle{definition}
\newtheorem{definition}{Definition}

\ifCLASSINFOpdf
\else
\fi

\usepackage{comment}
\usepackage{ifthen}
\newboolean{showcomments}
\setboolean{showcomments}{true}
\newcommand{\ychen}[1]{\ifthenelse{\boolean{showcomments}}
	{ \textcolor{red}{YC: #1}}}
\newcommand{\tongxin}[1]{\ifthenelse{\boolean{showcomments}}
	{ \textcolor{blue}{(#1)}}{}}

\usepackage{amsmath}
\allowdisplaybreaks[4]

\usepackage[ruled]{algorithm2e}

\begin{document}
	
	\title{Decentralized Provision of Renewable Predictions within a Virtual Power Plant}
	
	\author{Yue Chen,
		Tongxin Li,
		Changhong Zhao,
		and Wei Wei}% <-this % stops a space
%	\thanks{Y. Chen, W. Wei are with the State Key Laboratory of Power Systems, Department of Electrical Engineering, Tsinghua University, Beijing. (e-mail: cy11@tsinghua.org.cn, wei-wei04@mails.tsinghua.edu.cn)}
%	\thanks{T. Li is with Department of Computing and Mathematical Sciences, California Institute of Technology, Pasadena, CA. (email: {tongxin@caltech.edu})}
%	\thanks{C. Zhao is with the Department of Information Engineering, the Chinese University of Hong Kong, HKSAR, China. (email: chzhao@ie.cuhk.edu.hk)}}
	

	\markboth{Journal of \LaTeX\ Class Files,~Vol.~XX, No.~X, Feb.~2019}%
	{Shell \MakeLowercase{\textit{et al.}}: Bare Demo of IEEEtran.cls for IEEE Journals}



	\maketitle
	
	\begin{abstract}
The mushrooming of distributed energy resources turns end-users from passive price-takers to active market participants. To manage those massive proactive end-users efficiently, virtual power plant (VPP) as an innovative concept emerges. It can provide some necessary information to help consumers improve their profits and trade with the electricity market on behalf of them. One important information that is desired by the consumers is the prediction of renewable outputs inside this VPP.
Presently, most VPPs run in a centralized manner, which means the VPP predicts the outputs of all the renewable sources it manages and provides the predictions to every consumer who buys this information. We prove that by providing predictions, the social total surplus can be improved.
 However, when more consumers and renewables participate in the market, this centralized scheme needs extensive data communication and may jeopardize the privacy of individual stakeholders.
 In this paper, we propose a decentralized prediction provision algorithm in which consumers from each subregion only buy local predictions and exchange information with the VPP. Convergence is proved under a mild condition, and the demand gap between centralized and decentralized schemes is proved to have zero expectation and bounded variance. Illustrative examples show that the variance of this gap decreases with more consumers and higher uncertainty, and validate the proposed algorithm numerically. 

%The mushrooming of distributed renewable energy and storage endows consumers with higher flexibility in electricity markets. Traditionally, the market clearing is modeled via a Nash equilibrium problem, assuming consumers can predict the output of every renewable plant via information exchange, and then determine the market purchases which impact the clearing price. However, when more consumers join the market, global prediction needs extensive data communication and may jeopardize the privacy of individual players. In this paper, we propose a distributed prediction and decision (DPD) algorithm for a set of consumers in a local power market. Each consumer only predicts \textit{local} renewable output, and exchanges information with consumers in other microgrids/communities. Convergence is proved under a mild condition, and the gap between centralized (with global renewable predictions) and distributed methods (with local renewable predictions) is shown to have zero expectation with bounded variance. Illustrative examples validate the proposed algorithm numerically and show that the bound of gap variance decreases with more consumers and higher uncertainty.
	\end{abstract}
	

	% \fi
	\IEEEpeerreviewmaketitle
	
	\section*{Nomenclature}
	\addcontentsline{toc}{section}{Nomenclature}
	\subsection{Indices and Sets}
	\begin{IEEEdescription}[\IEEEusemathlabelsep\IEEEsetlabelwidth{${\underline P _{mn}}$,${\overline P _{mn}}$}]
		\item[$i, \mathcal{I} (\mathcal{I}_l)$] Index and set of consumers.
		\item[$l, \mathcal{L}$] Index and set of subregions.
		\item[$n, \mathcal{N}$] Index and set of generators.
		\item[$s, \mathcal{L}_g$] Index and set of lines.
		\item[$\mathcal{N}_l$] Set of consumers buying predictions on renewable output in subregion $l$.
		\item[$\mathcal{S}_i$] Set of predictable subregions for consumer $i$.
	\end{IEEEdescription}
	\subsection{Parameters}
	\begin{IEEEdescription}[\IEEEusemathlabelsep\IEEEsetlabelwidth{${\underline P _{mn}}$,${\overline P _{mn}}$}]
		\item[$\overline{w}, \sigma_{W_l}^2$] Expectation and variance of $W_l$.
		\item[$I, L$] Number of consumers and subregions.
		\item[$N, J$] Number of generators and fixed demand.
		\item[$c_n,d_n$] Cost coefficients of generator $n$.
		\item[$u,t$] Coefficients of consumer's utility function.
		\item[$D_j^f$] Fixed demand in the electricity market.
		\item[$\pi_{ns},\pi_{js},\pi_{vs}$]Line flow distribution factors.
		\item[$F_s$] Transmission capacity on line $s$.
		\item[$\alpha,\beta_0$] Coefficients of function indicating the relationship of electricity price and VPP demand, and we have $\gamma=I\alpha$, $\zeta=IL\alpha$.
		\item[$m$] Cost coefficient for certain precision level.
	\end{IEEEdescription}
		\subsection{Decision Variables}
		\begin{IEEEdescription}[\IEEEusemathlabelsep\IEEEsetlabelwidth{${\underline P _{mn}}$,${\overline P _{mn}}$}]
		\item[$D_i$] Consumption level of consumer $i$, and $D_i^{cen}, D_i^{dis}$ for centralized and decentralized schemes, respectively; $\overline{D}$ is the average value.
		\item[$g_n$] Power output of generator $n$.
		\item[$W_l$] Total renewable output in subregion $l$, a random variable and $w_l$ its realization.
		\item[$W_{il}^{pre}$] Prediction of $W_l$ for $i$, and $w_{il}^{pre}$ its realization.
		\item[$\lambda$] Electricity price for the virtual power plant.
		\item[$\epsilon_{il}$] Additive random noise in the prediction $W_{il}^{pre}$, and its variance is $\sigma_{\epsilon_{il}}^2$.
		\item[$\tau_{il}$] Prediction precision of consumer $i$ on renewable output in subregion $l$; and $\tau_{l}^{cen}, \tau_{l}^{dis}$ for centralized and decentralized schemes, respectively.
		\item[$\beta_l$] A signal in the decentralized algorithm.
		\end{IEEEdescription}

	\section{Introduction}

	\IEEEPARstart{T}{he} increasing penetration of distributed renewable energy has been considered as a promising solution to global warming and environmental pollution \cite{adefarati2016integration}. Despite the potential benefits, the volatility and intermittency of renewable energy also exerts challenges on power system operation \cite{phuangpornpitak2013opportunities}. Vast literature has been addressing this issue by stochastic \cite{xu2019data}, robust \cite{wang2016robust}, or distributionally robust \cite{wei2015distributionally} optimizations. When it comes to  distributed renewables, this problem becomes more complicated: Their capacities are usually small, so it is extremely hard for the operator to detect and prepare for their uncertainty beforehand \cite{liu2014system}. Facing these obstacles, there are three possible paradigms for designing energy markets 
to manage these distributed resources, i.e. peer-to-peer scheme, prosumer-to-grid scheme, and community-based scheme \cite{parag2016electricity}. Virtual power plant (VPP), whose initial version was the virtual utility in \cite{awerbuch2012virtual}, is a special case of the community-based approach, which integrates prosumers from different geographical locations. Therefore, it achieves a trade-off between  the peer-to-peer scheme and the prosumer-to-grid scheme, by being more organized than the former one while more flexible than the latter one which is restricted by the network structure. The main objective of VPP is to encourage DERs'/prosumers’ participation by improving their profits from the energy market \cite{zhang2018comprehensive}. Specially in this paper, the VPP operator can help prosumers gain higher profits in the electricity market by providing renewable output predictions. One practical example is the Enco Group in Netherlands and Belgium \cite{parag2016electricity}. %Despite the lack of universally accepted definitions, the essence of VPP is to promote system flexibility and economics by allowing different participants to virtually share their resources.
	A comprehensive review of VPP can be found in \cite{nosratabadi2017comprehensive}.
	
	In recent decades, VPP has captured great attention from both academia and industry. Basically, there are two ways for controlling devices and transmitting information within a VPP: the centralized approach and the decentralized approach \cite{zhang2018comprehensive}. Under the centralized scheme, there is a coordination centre in charge of all devices and participants. Information collection and transmission are conducted centrally, and global optimal strategies are made. The day-ahead optimal scheduling of VPP with energy storage, electrical and thermal energy resources, and demand response facilities were studied in \cite{zamani2016day}, and the point estimate method was used to depict the uncertainties. The imperialist competitive algorithm was proposed in \cite{kasaei2017optimal} to minimize the VPP's operating cost under renewable output, load demand, and market price uncertainties. The trading of VPP in both energy and reserve electricity markets were studied in \cite{baringo2018day} by stochastic adaptive robust optimization. The optimal automatic frequency restoration reserve for VPP was analyzed in \cite{camal2018optimal}. The robust capability curve indicating the allowable range of VPP was characterized in \cite{tan2020estimating}.

	 The centralized scheme can achieve global optimal solutions, but when it comes to large-scale distributed energy resources (DERs), it may become impracticable due to: 1) \emph{Computational Burden.} With an increasing number of DERs, the centralized scheme entails a tremendous amount of data acquisition and communication overhead.  2) \emph{Privacy}. The DERs are located in different subregions (SRs), owned by different stakeholders with conflicting economic interests, and they may be unwilling to provide private data for central control. Decentralized schemes can overcome these problems. A fully distributed ADMM-based algorithm for VPP problems was developed in \cite{chen2018fully}, and an application can be found in the demand response program of electric vehicles \cite{li2016admm}. Distributed primal-dual sub-gradient algorithm was proposed for coordinating DERs within a VPP via limited communication \cite{yang2013distributed}. Non-ideal communication network was considered in \cite{cao2017distributed}
	
	For both schemes, the gathering, processing and provision of information is critical as it is the foundation of the VPP operation. The interface between the VPP and DERs was designed with an information model in \cite{biegel2013information} so that the VPP can make full use of the aggregated flexibility of those resources. A general architecture as well as the Information and Communication Technology (ICT) infrastructure for VPP operation were presented in \cite{messinis2016ict}. A case was tested in \cite{kolenc2017performance} regarding the VPP communication system providing manual frequency restoration reserve service. Above works mainly focus on the physical construction of the VPP communication system, but the information to offer, the combination of information transmission and control schemes, and the resulting performance are also important issues that have not been formally studied. 
	This topic appears to be more important with the proliferation of distributed renewable energy. With large scale distributed renewable energy involved whose output is uncertain and volatile, the electricity price may fluctuate a lot making it hard for consumers to decide their optimal demand \cite{alahyari2020managing}. VPP can integrate prosumers from different subregions (SRs) and help them improve their profits from the electricity market \cite{zhang2018comprehensive} by services such as providing predictions of uncertain renewable outputs. We call this “\emph{uncertainty prediction provision}”.
	%The volatile output of renewable energy intensifies the uncertainty of electricity market prices, which will greatly influence the profit of the VPP and its participants \cite{alahyari2020managing}. 
	%Therefore, the VPP is desired to provide certain information (prediction) to help consumers decide their consumption levels.
	Providing such information incurs a cost related to information precision and thus leaves an optimal precision to be determined.
	
	This paper addresses the problem above by proposing a decentralized uncertainty prediction provision algorithm. The main contributions are two-fold: 
	
	1) \textbf{Model}. A two-level model is proposed to depict consumers' decision-making under distributed renewable uncertainties. In contrast to previous research, each consumer can purchase prediction information from the VPP to improve its knowledge of the uncertain renewable outputs via conditional expectation. The cost of this prediction information related to prediction precision is characterized.
	In the lower level, real-time electricity price is determined by a market-clearing problem with volatile renewables. The relationship between the price and the total demand inside VPP is further depicted as a piecewise linear curve. In the upper level, each consumer maximizes its conditional expected utility based on available predictions, and solves for an optimal prediction precision. We prove that providing predictions can boost social total surplus.
	
	2) \textbf{Algorithm}. A decentralized prediction provision (DPP) algorithm is presented. The region covered by the VPP is divided into several SRs, and the consumers inside each SR can only purchase local predictions of the renewable outputs in that SR. Meanwhile, it can exchange information with the VPP. The proposed decentralized scheme can overcome the drawbacks of centralized schemes such as extensive data exchange and possible violation of information privacy.
	We prove that the DPP algorithm converges under mild conditions, the gap of optimal strategies between centralized and decentralized cases has zero expectation and bounded variance. Simulation further shows that the variance of this gap decreases with more consumers and higher uncertainty.
	
	The model in this paper adopts a similar structure as \cite{wang2013privacy}. The differences are: This paper proposes a decentralized way to provide uncertainty predictions within a VPP, aiming at improving the profit of each consumer, while \cite{wang2013privacy} tries to minimize the cost of external energy imported to the system. Reference \cite{wang2013privacy} considers the privacy of consumers’ demand data, and assumes that the renewable generation in each cell is predictable, while this paper focuses on the privacy of renewable output data in each subregion taking into account the renewable uncertainty. Reference \cite{wang2013privacy} proposes a dual decomposition-based algorithm with decisions made at the cell level, while the proposed decentralized uncertainty provision algorithm derived from a  probabilistic perspective involving decision-making at the consumer level (corresponding to the participants/facilities level in \cite{wang2013privacy}).
	
In the rest of this paper, the market model and problem formulation are presented in Section \ref{sec:formulation}; the centralized benchmark for optimal power consumption is introduced in Section \ref{sec:centralized}; the distributed prediction provision algorithm is developed and compared with the centralized one in Section \ref{sec:distributed}; Numerical experiments are reported in Section \ref{sec:simulation}.
%\textbf{Notation}. We use $\mathbbm{E}[\cdot]$ to denote the expectation of random variables, and $\mathbbm{E}[\cdot|\cdot]$ the conditional expectation. We use $\mathbbm{D}(\cdot)$ to denote the variance of random variables, and $\mathrm{cov}(\cdot,.\cdot)$ the covariance. (CZ commented out for space.)

\section{Problem Formulation}\label{sec:formulation}
\subsection{System configuration and model setting}
\begin{figure}[t]
	\centering
	\includegraphics[width=0.9\columnwidth]{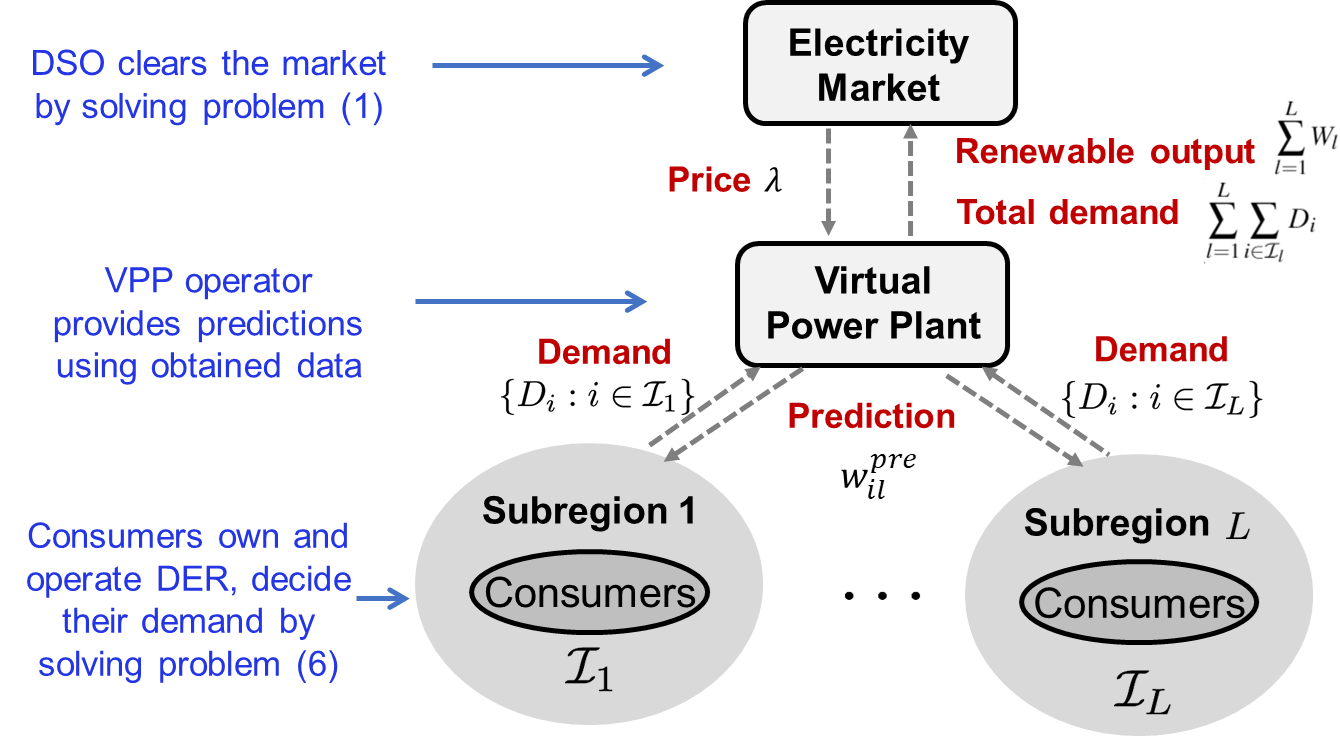}
	\caption{The hierarchical market model considered in this paper.}
    \vspace{-1mm}
	\label{fig:information flow}
\end{figure}

The problem considered in this paper involves three layers: distribution system operator (DSO), VPP operator, and consumers, as in Fig.\ref{fig:information flow}. The DSO clears the electricity market. The consumers own and operate the distributed energy resources (DERs). More accurate prediction of these DER outputs will help improve consumers’ profits and is desired. However, it could be impractical for each consumer to analyze the renewable output data and make prediction individually since it entails sophisticated techniques. Therefore, in our setting as in Fig.\ref{fig:information flow}, each subregion submits the output data of DERs in it to the VPP operator. The VPP operator makes predictions and provides them back to each subregion.
To be specific, first, the VPP operator provides predictions to consumers and consumers decide their optimal demand strategies. Then, the DSO collects data of the real-time renewable outputs and total demand from the VPP operator, together with data related to other demands and generators, and clears the electricity market centrally via problem \eqref{eq:market}.

%Consider a VPP participating in the distribution electricity market. 
The region covered by the VPP can be roughly divided into $L$ subregions indexed by $l \in \mathcal{L}:=\{1,...,L\}$ according to the location of distributed renewable sources. 
For each SR $l$, the aggregated volatile renewable output is uncertain and can be represented as a  random variable $W_l$ in $\mathbbm{R}^+$ with unknown distribution.  Denote  by $\overline{w}_l:=\mathbbm{E}(W_l)$ its expectation and  $\sigma_{W_l}^2:=\mathbbm{D}(W_l)$ its variance. Each SR consists of a set of consumers indexed by $i\in \mathcal{I}_l:=\{1,\ldots,I_{l}\}$, and the set of all consumers inside the VPP is $\mathcal{I}=\cup_{l=1}^L \mathcal{I}_l$. Every consumer can adjust its consumption level to maximize its (expected) net utility, and the demand of consumer $i$ in SR $l$ is $D_i,i\in \mathcal{I}_l$\footnote{The consumer's demand is determined by the consumer's decision-making problem \eqref{eq:consumer} and we will give a detailed description in Section~\ref{sec:decision}.}.
The VPP gathers the demand requirements from all consumers, which is $\sum_{l=1}^L \sum_{i \in \mathcal{I}_l} D_i$, and purchases electricity from the distribution power market. It is worth noting that, in contrast to traditional retailer, the VPP does not make profit directly from selling electricity but acts as a coordinator and provides supporting information/prediction.

\subsection{Market clearing model}
The prevalence of price-sensitive flexible loads and volatile output of renewable energy intensifies the uncertainty of electricity market prices, which will greatly influence the profit of the consumers in the VPP. In this context, consumer’s behavior as well as the real-time market stability has become a crucial topic needs further investigation \cite{zhang2018comprehensive}. Here, we consider the hour-head bidding problem as \cite{zhao2014convergence,fang2015coupon,pei2016optimal}, and draw some in-depth insights from a theoretical point of view.
 
Suppose there are $N$ thermal generators indexed by $n \in \mathcal{N}=\{1,...,N\}$ in the distribution electricity market. The output of the $n$-th generator is $g_n$, and its cost function is $f_n(g_n)=c_ng_n^2+d_ng_n$, where $c_n$ and $d_n$ are non-negative cost coefficients. The cost of a renewable unit is usually close to zero, and therefore is neglected here. The market minimizes the total generation cost in \eqref{eq:market.1} subject to power balancing constraint \eqref{eq:market.2}, capacity limits \eqref{eq:market.3} and network constraints \eqref{eq:market.4}:
\begin{subequations}
	\label{eq:market}
	\begin{align}
	\label{eq:market.1}
	\mathop{\min}_{g_1,\ldots,g_{N}} & \sum \limits_{n=1}^{N} (c_ng_n^2+d_ng_n) \\
	\text{s.t.} & \sum \limits_{n=1}^{N} g_n + \sum \limits_{l=1}^{L} W_l = \sum \limits_{l=1}^{L} \sum \limits_{i\in \mathcal{I}_l} D_{i} + \sum \limits_{j=1}^J D_j^f: \lambda\label{eq:market.2} \\
	~ & \underline{g}_n \le g_n \le \overline{g}_n, \forall n \in \mathcal{N} \label{eq:market.3} \\
	~ & \Big| \sum \limits_{n=1}^N \pi_{ns}g_n- \sum \limits_{j=1}^J \pi_{js}D_j^f + \pi_{vs} \left(\sum \limits_{l=1}^{L} W_l - \sum \limits_{l=1}^{L} \sum \limits_{i\in \mathcal{I}_l} D_{i}\right) \Big|   \nonumber\\
	~ & \le F_s,\forall s \in \mathcal{L}_g \label{eq:market.4} 
	\end{align}
\end{subequations}
Here, $D_j^f,\forall j=1,...,J$ represents the equivalent fixed demand other than the purchase of VPP in the distribution electricity market. Note that from the perspective of VPP, different from the random renewable output, the inaccurate load forecast of $D_j^f$ is an external uncertainty, reflected in the term $\beta_0$. In general, the accuracy of hourly ahead load forecast is much higher than that of renewable prediction \cite{zhang2015day}, so we assume $D_j^f,\forall j=1,...,J$ are constants and only consider the uncertainty of renewable generation.
$\lambda$ comes from the dual variable of the power balancing condition \eqref{eq:market.2}, giving the electricity price for the VPP (also for every consumer inside this VPP).

To facilitate further analysis, we first characterize the relationship between price $\lambda$ and the VPP's total demand $\sum_{l=1}^L \sum_{i \in \mathcal{I}_l} D_i$ as a piecewise linear curve, whose particular segment around the market clearing point is:
%\begin{align*}
%2c_n g_n + d_n - \lambda &= 0, \ \text{for all } n=1,\ldots,N, \\
%\sum_{n=1}^{N} g_n + \sum_{l=1}^{L} W_l &= %\sum_{l=1}^{L} \sum_{i \in \mathcal{I}_l} D_{i},
%\end{align*}
\bq
\label{eq:marketprice}
\lambda(D,W) = \alpha \left(\sum_{l=1}^{L} \sum_{i \in \mathcal{I}_l} D_{i} - \sum_{l=1}^{L} W_l\right)+\beta_0
\eq
with constant coefficients $\alpha>0$, $\beta_0$. 

\textbf{Remark:} 
(i) Price $\lambda$ as the Lagrangian multiplier of constraint \eqref{eq:market.2} can be solved from the KKT condition at the primal-dual optimum of \eqref{eq:market}. It can be proved under mild conditions (e.g., at optimum, at least one generator does not reach its capacity limit and is not causing congestion) that the dual-optimal $\lambda$ is unique (for which we skip the proof) and is a piecewise linear function of $\left(\sum_{l=1}^{L} \sum_{i \in \mathcal{I}_l} D_{i} - \sum_{l=1}^{L} W_l\right)$ \cite[p. 95]{still2018lectures}. Such a piecewise linear structure is also observed in power market literature, e.g., \cite{zhao2014convergence}.
(ii) Since one VPP only accounts for a small fraction of the grid's total load, it is reasonable to assume that \eqref{eq:market} has a fixed set of binding constraints at optimum, despite the varying renewable generation and demand of the VPP. Therefore, $\lambda$ stays in the same segment of the piecewise linear curve, whose coefficients $\alpha$, $\beta_0$ can be determined a priori from historical data.
(iii) We assume that the total demand $\sum_{l=1}^L \sum_{i \in \mathcal{I}_l} D_i$ is sufficiently higher than the total renewable output $\sum_{l=1}^{L} W_l$, which keeps $\lambda(D,W)$ positive.

\subsection{Prediction information of renewable sources}

As mentioned above, the renewable outputs $W_1,\ldots,W_{L}$ are random variables, and thus, the electricity price $\lambda(D,W)$ is also volatile and uncertain. The VPP may try to offer prediction of the output of renewable units $W_1,\ldots,W_{L}$ to help consumers enhance their profits. Consumers need to pay for this information, and $\mathcal{N}_l$ is the set of consumers buying predictions on renewable output in SR $l$.  If consumer $i$ buys the information for SR $l$, then it will get a forecast of $W_l$ from $W_{il}^{pre}=W_l+\varepsilon_{il}$, where $\varepsilon_{il}$ is additive random noise. $W_{il}^{pre}$ is also a random variable, and $w_{il}^{pre}$ is the actual forecast consumer $i \in \mathcal{N}_l$ gets. We assume:
\begin{enumerate}[label=A\arabic*:]
	\item  Each SR has the same number of identical consumers, i.e., $I_l=I$ for all $l=1,\ldots,L$.
	\item $\{\varepsilon_{il}:i\in\mathcal{N}_{l}\}$ are independent and identically distributed (i.i.d.), each $\varepsilon_{il}$ has zero expectation and is independent of $W_l$. Denote by $\sigma_{\varepsilon_{il}}^2$ the variance of noise $\varepsilon_{il}$.
\end{enumerate}
%\begin{enumerate}[label=A\arabic*:]
%\item  Each microgrid has the same number of identical consumers, i.e., $I_l=I,\forall l$.
%\item  $\mathbbm{E}[W_l|w_{il}^{pre}]$ is affine in $w_{il}^{pre}$.
%\end{enumerate}
Although different SRs can have different numbers of consumers and the consumers inside each SR may also be heterogeneous, we make assumption A1 to derive a theoretically sound result. In practice, we can adjust the region division so that every SR has an equal number of consumers as well as a representative to model all the consumers inside it.

Consumer $i$ utilizes $w_{il}^{pre}$ to improve its knowledge about $W_l$.
It is well known that conditional expectation $\mathbbm{E}[W_l|W_{il}^{pre}=w_{il}^{pre}]$ is the best estimator of $W_l$ \cite[p. 346]{grimmett2001probability}, which, however, may be a complicated nonlinear function of $w_{il}^{pre}$ if $(W_l, W_{il}^{pre})$ is not subject to bivariate normal distribution \cite{tntech}. Therefore, in this paper we adopt the \emph{best linear estimator} of $W_l$ as\footnote{We still use $\mathbbm{E}$ to denote the best linear estimator although it can be an approximate of the actual expectation.}
\begin{eqnarray}
\mathbbm{E}[W_l|W_{il}^{pre}=w_{il}^{pre}] &=& A_{1i,l}+A_{2i,l} w_{il}^{pre} \nonumber
\end{eqnarray}
where $A_{1i,l}$ and $A_{2i,l}$ are given in the following lemma.
\begin{lemma}
	\label{lemma1}
	\begin{align}
	A_{2i,l} = \frac{\sigma_{W_l}^2}{\sigma_{W_l}^2+\sigma_{\varepsilon_{il}}^2}, \;\; A_{1i,l}=(1-A_{2i,l})\overline{w}_l
	\end{align}
	%\vspace{0.1mm}
\end{lemma}

The proof of Lemma \ref{lemma1} can be found in Appendix.\ref{apen-lemma1}, and it is indeed consistent with known results about the best linear estimator, e.g., \cite{tntech}.
When $\sigma_{\varepsilon_{il}}^2 \to \infty$, we have $\mathbbm{E}[W_l|W_{il}^{pre}=w_{il}^{pre}]=\overline{w}_l$. It means that the VPP's forecasting technique is poor and the prediction $w_{il}^{pre}$ offers no information, so the consumer simply uses the expected renewable output. 
With higher penetration of renewables, forecasting plays an increasingly important role in enhancing consumers' profits.
 However, improving the precision of the forecasts from VPP may generate costs, which will be split among the consumers.
 To depict this cost, we borrow the concept of prediction precision in economics \cite{vives2010information} defined as follows.
with the related cost given in \eqref{eq:precision-cost}.

\begin{definition} (Prediction precision~\cite{vives2010information}) The coefficient $A_{2i,l}$ in Lemma \ref{lemma1} is formally defined as \textit{prediction precision} of consumer $i$ on renewable output in SR $l$:
	\begin{align}
	\label{eq:precison}
	\tau_{il}:=\frac{\sigma_{W_l}^2}{\sigma_{W_l}^2+\sigma_{\varepsilon_{il}}^2}
	\end{align}
\end{definition}
The smaller $\sigma_{\varepsilon_{il}}^2$ is, the more accurate the prediction is. The concept $\tau_{il}$ given in \eqref{eq:precison} decreases with $\sigma_{\varepsilon_{il}}^2$: when $\sigma_{\varepsilon_{il}}^2 \to 0$, $\tau_{il} \to 1$ and we have $W_{il}^{pre}=W_l$ which is the most accurate case; when  $\sigma_{\varepsilon_{il}}^2 \to \infty$, $\tau_{il} \to 0$ and we have $W_{il}^{pre} \approx \varepsilon_{il}$ which is the most inaccurate case. $\tau_{il}$ varies in $[0,1]$ and a larger $\tau_{il}$ indicates a higher prediction precision. Improving prediction technique may incur costs \cite{vives2010information}, which is paid by all consumers in $\mathcal{N}_l$ and given by
\begin{align}
\label{eq:precision-cost}
h_l(\tau_{il}):=\frac{\check{m}}{\sigma_{W_l}^2}\frac{\tau_{il}}{1-\tau_{il}}
\end{align}
where $\check{m}$ is the cost coefficient for achieving certain precision level, and $\check{m}:=m/|\mathcal{N}_l|$. Denote by $\tau_l := {\sum_{i \in \mathcal{N}_l} \tau_{il}}/{|\mathcal{N}_l|}$. The prediction cost of renewable output in SR $l$ for each consumer $i$ is related to three factors, the number of consumers buying this information $|\mathcal{N}_l|$, the prediction accuracy $\tau_{il}$, and the level of renewable uncertainty $\sigma_{W_l}^2$. The term $\check{m}$ decreases with $|\mathcal{N}_l|$, indicating that with more consumers involved, the cost per consumer goes down. The term $\tau_{il}/(1-\tau_{il})$ increases with $\tau_{il}$, meaning that the consumer needs to pay more for a higher prediction precision. A higher $\sigma_{W_l}^2$ would make it harder to predict the renewable output accurately so the cost increases.

\subsection{Decision making of consumers}
\label{sec:decision}

The demand of consumer $i$ is $D_i$, and it pays at the electricity price $\lambda(D,W)$ in~\eqref{eq:marketprice}. The real-time price $\lambda(D,W)$ depends not only on the joint consumption of all consumers, but also on the realization of random outputs $W_1,\ldots,W_{L}$ which is not known to the consumers. Therefore, consumers need to predict future renewable outputs for making their decisions. The consumer can buy these predicitons from the VPP and we denote by $\mathcal{S}_{i}$ the set of \textit{predictable subregions} whose predictions are available to the consumer $i$. 
In other words, the $i$-th consumer knows the predictions $\{w_{il}^{pre}: l\in\mathcal{S}_i\}$ before determining $D_{i}$, the amount of electricity to buy. 
Define an event $\mathcal{A}_i:=\{W_{il}^{pre}=w_{il}^{pre}, \  l \in \mathcal{S}_{i}\}$ for all consumers $i\in\mathcal{I}$. 

The utility function of consumer $i$ is $U_{i}(D_{i})=-uD_{i}^2/2+tD_{i}$ where $u,t>0$ are constants. The consumer $i$ chooses its demand level to maximize its net utility, i.e., the utility minus the electricity purchase cost. Let the conditional expectation  $\mathbbm{E}[\lambda(D,W)|\mathcal{A}_i]$ be the prediction of price for each consumer $i\in\mathcal{I}$. This leads to the following decision rule of $D_i, i \in \mathcal{I}$:
% However, since $\lambda(D,W)$ is unknown, consumer $i$ can only estimate its value based on the predictions, as a conditional expectation $\mathbbm{E}[\lambda(D,W)|\mathcal{A}]$.
% Here, $\mathcal{S}_{i}, \forall i\in\mathcal{I}$ is the set of microgrids consumer $i$ has predictions on, which means if $i \in \mathcal{N}_l$ then $l \in \mathcal{S}_i$. It is worth noting that, $\mathcal{S}_{i}$ does not necessarily contain predictions over all renewable outputs $W_l,\forall l$.
\begin{align}
\label{eq:consumer}
\mathop{\max}_{D_i} ~ \pi_i:= U_i(D_{i}) &
-D_{i}\mathbbm{E}[\lambda(D,W)|\mathcal{A}_i]
\end{align}

We consider the case where there is such a large number of consumers in the market that the impact of an individual consumer's strategy on the electricity price $\lambda(D,w)$ is negligible. It means that $\mathbbm{E}[\lambda(D,w)|\mathcal{A}_i]$ can be treated as a given constant in problem \eqref{eq:consumer} \cite{faqiry2016double}. Thus, letting the derivative be zero, the $i$-th consumer's decision $D_{i}^*$ satisfies:
\begin{align}
\label{eq:condition}
t-uD_{i}^* = \mathbbm{E}[\lambda(D,W)|\mathcal{A}_i]
\end{align}
%	\vspace{-3mm}
\subsection{Predictable subregions: centralized vs decentralized}
\label{sec: two_models}

\begin{figure}[t]
	\centering
	\includegraphics[width=0.85\columnwidth]{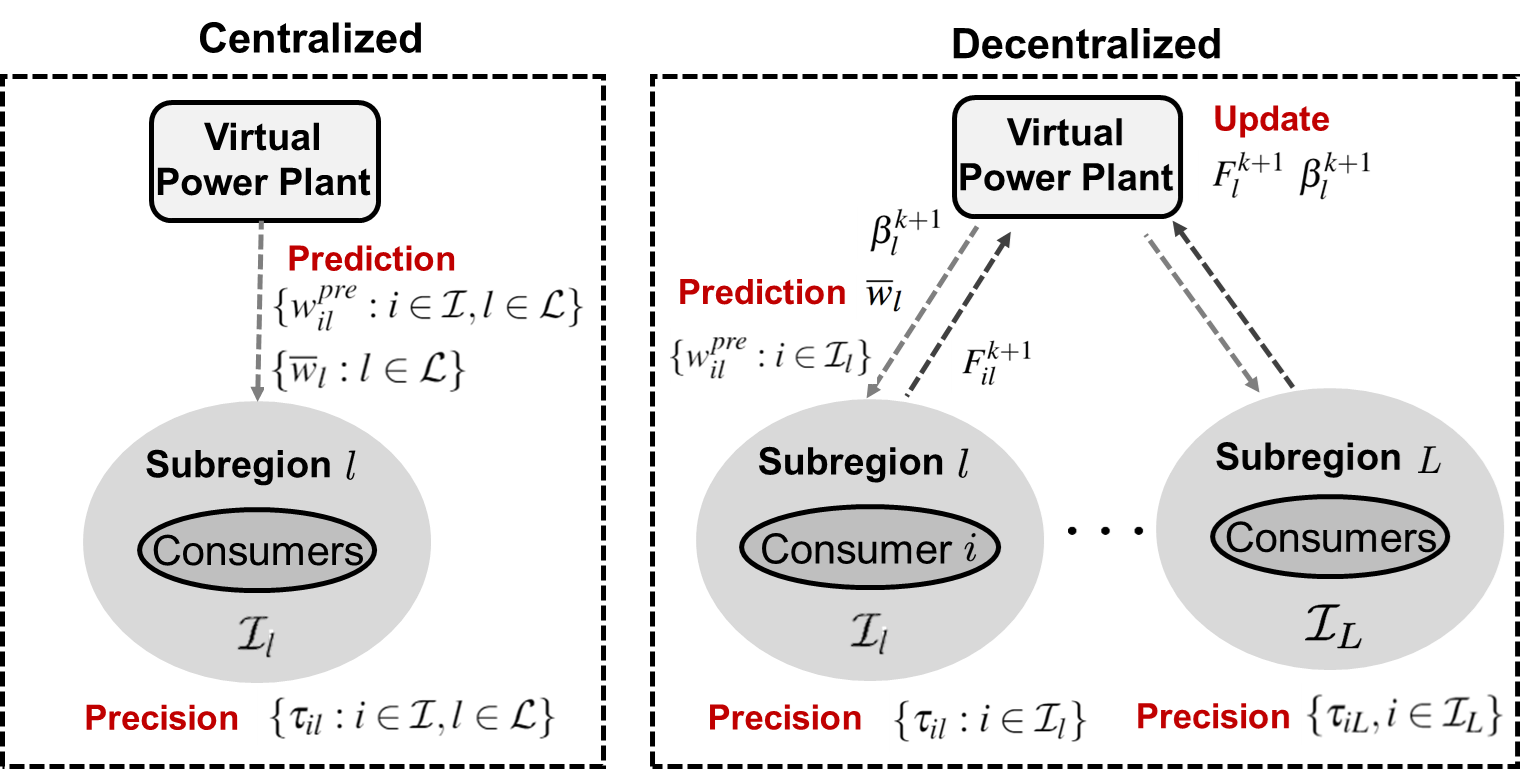}
	\caption{Centralized (left) versus decentralized (right) prediction provision.}
	\vspace{-1mm}
	\label{fig:comparison}
\end{figure}

The choice of $\{D_i^*:i\in\mathcal{I}\}$ will be influenced by the sets of SRs $\{\mathcal{S}_i:i\in\mathcal{I}\}$ that the consumers can get predictions on. In the following, we compare two VPP's information provision ways: \emph{centralized} and \emph{decentralized}. Under the centralized scheme, each consumer is able to know the predictions of all the SRs from the VPP. Under the decentralized scheme, each consumer can only get local predictions from the SR it belongs to but can exchange some information with other SRs via the VPP. The centralized scheme serves as a benchmark in our analysis, which is presented in the sequel.

\textbf{\textit{Centralized scheme.}}
 Suppose \emph{every} consumer get predictions of uncertain renewable outputs in all the SRs, as shown in Fig.\ref{fig:comparison} (left). Therefore, the predictable set $\mathcal{S}_{i}=\mathcal{L}$ for all $i\in \mathcal{I}$ and $\mathcal{N}_l=\mathcal{I}$ for all $l \in \mathcal{L}$. Specifically, each consumer $i\in\mathcal{I}$ has predictions $\{w_{il}^{pre}:l\in\mathcal{L}\}$ and expectations $\{\overline{w}_{l}: l\in\mathcal{L}\}$. This centralized scheme in Section \ref{sec:centralized} serves as a benchmark for the decentralized scheme in Section \ref{sec:distributed}.

\textbf{\textit{Decentralized scheme.}}
The centralized benchmark above could be hard to implement  because of: 1) \emph{Communication burden}. The centralized scheme requires the information transmission of all SRs to all consumers, which will incurs huge challenge on the data acquisition and communication infrastructure. 2) \emph{Information privacy}. Different SR may want to keep its data privacy, and are not willing to share all those information with other SRs.
Therefore, we consider a more practical case where only local predictions are available, as shown in Fig.\ref{fig:comparison} (right). Specifically, the set of predictable SRs is $\mathcal{S}_i=\{l\}$ for every consumer $i\in \mathcal{I}_l$. Each consumer $i\in\mathcal{I}_l$ \textit{only} has prediction $\{w_{il}^{pre}\}$ and expectation $\overline{w}_{l}$ of the subregion $l$ where it is located.

We next present a centralized model with optimal affine decision policy based consumers as a benchmark.

\section{Centralized Scheme with an Affine Policy}\label{sec:centralized}
In this section, we analyze consumers' strategies under the centralized scheme as well as the optimal prediction precision.
In principle, consumers' strategy can be an arbitrary function of their predictions. Here, we restrict  the scope of this study to the affine policies, which is commonly used in uncertain related studies \cite{bertsimas2010optimality,lorca2016multistage}. Specifically, consider the following strategy of consumer $i \in \mathcal{I}$:
\begin{align}
\label{eq:consumer-stratetgy}
D_{i}^{cen*}:=b_1+\sum \limits_{l=1}^{L} b_{2i}^{l}(w_{il}^{pre}-\overline{w}_l)
\end{align}
where $b_1$ and $\{b_{2i}^{l}:i\in\mathcal{I},l\in\mathcal{L}\}$ are constants to be determined. 
Note that $\{\epsilon_{il}:i\in\mathcal{I}\}$ are i.i.d. with $\mathbbm{E}[b_{2i}^l\epsilon_{il}]=0$ and assume variance $\mathbbm{D}(b_{2i}^l\epsilon_{il})$ is uniformly bounded for all $i$ and $l$.  Consider a big group $\mathcal{N}_l$ of consumers so that the law of large numbers can be applied to obtain their average demand:
\begin{align*}
\overline{D}^{cen*}
:=& \sum \limits_{l=1}^L \frac{1}{|\mathcal{N}_l|} \sum_{i \in \mathcal{N}_l} b_{2i}^{l} (w_{il}^{pre}-\overline{w}_l)  +b_1\\
= & \sum \limits_{l=1}^L \frac{\sum_{i \in \mathcal{N}_l} b_{2i}^{l}}{|\mathcal{N}_l|} (w_l-\overline{w}_l)+ b_1 + \sum \limits_{l=1}^L \frac{\sum_{i \in \mathcal{N}_l} b_{2i}^{l} \varepsilon_{il}}{|\mathcal{N}_l|}\\
=& b_1+\sum \limits_{l=1}^L b_{2}^{l}(w_{l}-\overline{w}_l)
\end{align*}
where $b_2^{l}:={\sum_{i \in \mathcal{N}_l} b_{2i}^{l}}/|\mathcal{N}_l|={\sum_{i \in \mathcal{N}_l} b_{2i}^{l}}/{(IL)}$. Then, we can get every consumer's optimal strategy under given prediction precision by giving the following proposition.
\begin{proposition}
	\label{prop1}
	Given prediction precision $\tau_{il},\forall l \in \mathcal{L}, i \in \mathcal{N}_l$, consumer $i \in \mathcal{I}$'s optimal demand under centralized manner is given by \eqref{eq:consumer-stratetgy} with
	\begin{align}
	b_{2i}^{l} =   ~& \frac{\alpha \tau_{il}}{\zeta \tau_{l} + u}, \forall l \in \mathcal{L}, \nonumber\\
	b_1 =  ~ & \frac{t-\beta_0+\sum \limits_{l=1}^L \alpha \overline{w}_l}{\zeta+u}
	\end{align}
	where $\zeta=IL\alpha$.
%	where $\alpha={1}/{\sum \limits_{n=1}^N \frac{1}{2c_n}}$, $\beta_0={\sum \limits_{n=1}^N \frac{d_n}{2c_n}}/{\sum \limits_{n=1}^N \frac{1}{2c_n}}$, $\zeta=IL\alpha$.
	%\vspace{1mm}
\end{proposition}

The proof of Proposition \ref{prop1} is in Appendix.\ref{apen-prop1}. It gives the optimal strategy of a consumer under certain prediction precision. Together with the optimal condition \eqref{eq:condition}, consumer $i \in \mathcal{I}$'s expected net utility is
\begin{align}
\mathbbm{E}[\pi_{i}]:= ~& \mathbbm{E}[\frac{1}{2}u(D_{i}^{cen*})^2] \nonumber\\
= ~ & \frac{u}{2}\left[(\frac{t-\beta_0+\alpha \sum \limits_{l=1}^L \overline{w}_l}{\zeta+u})^2+\sum \limits_{l=1}^L \frac{\alpha^2\sigma_{W_l}^2\tau_{il}}{(\zeta \tau_{l}+u)^2}\right]
\end{align}
Note that $\mathbbm{E}[\pi_{i}]$ increases with prediction precision $\tau_{il}$ of consumer $i$, but decreases with the average precision $\tau_l$ across $\mathcal{N}_l$. This aligns with the intuition that a consumer earns more by having more information about renewables and earns less if all of its competitors know more information than itself. The optimal prediction precision of a consumer is given by Proposition \ref{prop2}, whose proof is deferred to Appendix.\ref{apen-prop2}.

\begin{proposition}
	\label{prop2}
Under the centralized scheme, the optimal prediction precision for the renewable source in subregion $l$ is
	\begin{align}
	\label{eq:prediction-accuracy-cen}
	\tau_{l}^{cen*}=\mathop{\max}\left(0,\frac{ \sigma_{W_l}^2-\sqrt{2mu/(IL)}/\alpha}{\sigma_{W_l}^2+ \sqrt{2mIL/u}}\right)
	\end{align}
	%\vspace{0.1mm}
for \emph{every} consumer $i \in \mathcal I$.
\end{proposition}

Following Proposition \ref{prop2}, we obtain the optimal average demand over all the consumers in $\mathcal{N}_l,\forall l \in \mathcal{L}$ as:
\begin{align}
\label{eq:opt-d-cen}
D_l^{cen*} = b_1 + \sum \limits_{j=1}^L \frac{\alpha \tau_{j}^{cen*}}{\zeta \tau_{j}^{cen*} + u}(w_j-\overline{w}_j)
\end{align}
For renewable uncertainty (variance) $\sigma_{W_l}^2$ so small that $\tau_{l}^{cen*}=0$, a consumer tends to make decisions just based on expectation $\overline{w}_l$. In contrast, as $\sigma_{W_l}^2 \to \infty$, a consumer's best option is $\tau_{l}^{cen*} \to 1$, i.e., to make as accurate prediction as possible.

Moving further, we discuss the market impact of providing predictions. First, the concept of total surplus is introduced to quantify market efficiency.
\begin{definition}
Given the renewable output $w_l,\forall l \in \mathcal{L}$ and each consumer's demand $D_i,\forall i \in \mathcal{I}_l, l \in \mathcal{L}$, the \emph{total surplus} of the market is the sum of consumer surplus and producer surplus:
\begin{align}
\label{eq:total-surplus}
    TS(D,W)=\sum \limits_{l=1}^L \sum \limits_{i \in \mathcal{I}_l} U_i(D_i)-\int_{0}^{\sum \limits_{l=1}^L \sum \limits_{i \in \mathcal{I}_l}D_i} \lambda(x,W)dx
\end{align}
\end{definition}

If all consumers know exactly the renewable outputs $w_l,\forall l \in \mathcal{L}$, denote their optimal demands as $D_i^0,\forall i \in \mathcal{I}$, and the corresponding total surplus as $TS(D^0,W)$. If no effort is devoted to improving the prediction techniques and consumers simply use the expected renewable output $\overline{w}_l,\forall l \in \mathcal{L}$, denote their optimal demand as $D_i^{1},\forall i \in \mathcal{I}$, and the corresponding total surplus as $TS(D^1,W)$. Given a renewable output profile $w _l,\forall l \in \mathcal{L}$, the conditionally expected total surplus is $\mathbbm{E}[TS(D^{cen*},W)|W=w]$. 
Theorem \ref{prop7} compares market efficiency of these three cases.
\begin{theorem}
\label{prop7}
If assumptions A1, A2 hold, we have
\begin{align}
    \mathbbm{E}[TS(D^0,W)] \ge \mathbbm{E}[\mathbbm{E}[TS(D^{cen*},W)|W=w]] \ge \mathbbm{E}[TS(D^1,W)] \nonumber
\end{align}
\end{theorem}
The proof of Theorem \ref{prop7} is in Appendix.\ref{apen-prop7}. With complete information about the renewable output, $\mathbbm{E}[TS(D^0,W)]$ represents the optimal total surplus. However, when we cannot know the renewable outputs accurately, a total surplus loss is incurred. By providing improved predictions using the centralized scheme, this loss can be reduced.

\section{Decentralized Scheme with an Affine Policy}\label{sec:distributed}

\subsection{Mechanism description}
Due to the reasons mentioned in Section~\ref{sec: two_models}, \emph{Decentralized Prediction Provision} (DPP) is highly desired, in which every consumer can only get access to the predictions of local renewable output of the SR it lies in, i.e. $\mathcal{S}_{i}=\{l\},\forall i \in \mathcal{I}_l$. In this paper, we propose a DPP mechanism in \textbf{Algorithm 1}, whose information flow is outlined in  Fig.\ref{fig:comparison} (right). Take the consumers in SR $l$ for example:

\textbf{Step 1:} (Initialization) The VPP generates the predictions $w_{il}^{pre}$ for all $i \in \mathcal{I}_l$ given the precision set by \eqref{eq:optimal-dis-pa}, and set $\beta_l=0$ for all SR $l$ (the intuition of $\beta_l$ will be explained later).

\textbf{Step 2:} Every consumer in $ \mathcal{I}_l$ receives a $\beta_l$ and local prediction information $w_{il}^{pre}$ from the VPP.

\textbf{Step 3:} With  $\beta_l$ and $w_{il}^{pre}$, it accordingly decides its optimal demand $D_{il}^{dis*}$ based on \eqref{eq:individualaffineDPD} and \eqref{eq:dis-consumer}\footnote{For the decentralized case, the optimal strategy and the demand of consumers need be differentiated across SRs.}. It also reports $F_{il}:=D_{il}^{dis*}-w_{il}^{pre}/I$ back to the VPP.

\textbf{Step 4:} The VPP sums up all $F_{il}$ from consumers $i \in \mathcal{I}_l$ and gets the result $F_l$ for SR $l$. Then it further updates $\beta_l$ according to $F_j$ from all other SR $j \in \mathcal{L}, j\neq l$, and sends the updated $\beta_l$ back to consumers in SR $l$. This process repeats imitatively until convergence.

%The microgrid operator receives $\beta_l$ from an information platform and broadcasts it to every consumer in $ \mathcal{I}_l$. The consumer $i\in\mathcal{I}_l$ also makes a prediction $w_{il}^{pre}$ about the renewable output of the local microgrid $l$. It accordingly decides its optimal demand $D_{il}^{dis*}$ via \eqref{eq:individualaffineDPD} and \eqref{eq:dis-consumer}.\footnote{For the distributed case, the optimal strategy and the demand of consumers need be differentiated across microgrids. Therefore, we have the microgrid index $l$ in the subscript of $D_{il}^{dis*}$.} Then the consumer $i$ reports $F_{il}:=D_{il}^{dis*}-w_{il}^{pre}/I$ to the microgrid $l$. The microgrid operator sums up all $F_{il},\forall i \in \mathcal{I}_l$ and returns the result $F_l$ to the information platform. The platform updates $\beta_l$ according to $F_j$ received from all other microgrids $j \in \mathcal{L}, j\neq l$ and sends updated $\beta_l$ back to the microgrid $l$ and this process repeats imitatively until convergence.

\begin{algorithm}[t]
\small
	\caption{\textbf{D}ecentralized \textbf{P}rediction \textbf{P}rovision} 
	%\LinesNumbered %要求显示行号
	\KwIn{input parameters $c_n, d_n,\forall n \in \mathcal{N}$, $I$, $L$, $\varepsilon$.}%输入参数
	\KwOut{optimal demand ${D}_{il}^{dis*},\forall l \in \mathcal{L}, \forall i \in \mathcal{I}_l$.}%输出
	\textbf{Initialization:} $\beta_l=0,\forall l \in \mathcal{L}$, $k=0$\; %\;用于换行
	\Repeat{$\mathop{\max}_{l \in \mathcal{L}}|F_l^{k+1}-F_l^k| \le \varepsilon$}{
		iteration $k++$
		
		\textbf{VPP:}
		
		\For{$l=1;l \le L$} 
		{
		\textbf{consumers in SR $l$:}
		
			\For{$i=1; i \le I$}
			{
				consumer $i$ gets $\beta_l^k$ and a predict $w_{il}^{pre}$.
				
				update $D_{il}^{k+1}, \forall i \in \mathcal{I}_l$ by \eqref{eq:individualaffineDPD}~\eqref{eq:dis-consumer}, and
				$$F_{il}^{k+1}:=D_{il}^{k+1}-w_{il}^{pre}/I,\forall i \in \mathcal{I}_l$$}
			$F_l^{k+1} := \sum_{i \in \mathcal{I}_l} F_{il}^{k+1} $
		}
		update $\beta_l^{k+1}$ by
		$$
		\beta_l^{k+1} := \beta_0 + \alpha \sum_{j \ne l} F_j^{k+1} , \forall l \in \mathcal{L}
		$$
	}
\end{algorithm}

\subsection{Decision in each iteration}
For SR $l$, the relationship between its total demand $\sum_{i \in \mathcal{I}_l} D_{il}$ and the market price $\lambda$ in \eqref{eq:marketprice} can be represented as below. The information related with other SRs are compressed in $\beta_l$.
\begin{align}
\lambda(D_l,W_l)\!=\!\alpha \left(\sum_{i \in \mathcal{I}_l} D_{il}-W_l \right) \!+\! \underbrace{\alpha \!\! \!\sum \limits_{j \ne l, j \in \mathcal{L}} \! \left(\sum \limits_{i \in \mathcal{I}_j} D_{ij} -W_{j}\right)+\beta_0}_{\beta_l} \nonumber
\end{align}
 
Ideally $\beta_l$ should contain $W_j$ for all $j \in \mathcal{L}, j \neq {l}$, which however are random and it is impossible to know their exact value $w_j$ in advance. In \textbf{Algorithm 1} we use $\sum_{i \in \mathcal{I}_j} w_{ij}^{pre} / I$ to approximate $w_j$ and calculate $\beta_l$. When $I$ is large, this can be a good approximation.

\textbf{Remark}: Note that in practice, there might be some customers without consumption flexibility. In that case, denote the set of fixed demand in SR $l \in \mathcal{L}$ as $\mathcal{I}_l^f$, then the market price   can be represented as
\begin{align}
\lambda(D_l,W_l)=~ & \alpha \left(\sum_{i \in \mathcal{I}_l/\mathcal{I}_l^f} D_{il}-W_l \right) \nonumber\\
~ & +\underbrace{ \alpha \sum \limits_{i \in \mathcal{I}_l^f} D_{il}+\underbrace{\alpha \!\! \!\sum \limits_{j \ne l, j \in \mathcal{L}} \! \left(\sum \limits_{i \in \mathcal{I}_j} D_{ij} -W_{j}\right)+\beta_0}_{\beta_l}}_{\hat \beta_l} \nonumber
\end{align}
After SR $l \in \mathcal{L}$ gets $\beta_l^k$, it first let $\hat \beta_l=\beta_l+\alpha \sum \nolimits_{i \in \mathcal{I}_f} D_{il}$ and $\hat I=I-|\mathcal{I}_l^f|$, then the flexible load $i \in \mathcal{I}_l/\mathcal{I}_l^f$  updates $D_{il}^k$ according to \eqref{eq:individualaffineDPD}~\eqref{eq:dis-consumer} using $\hat \beta_l^k$ to replace $\beta_l$, $\hat I$ to replace $I$. Finally, all consumers $i \in \mathcal{I}$ calculate $F_{il}^{k+1}=D_{il}^{k+1}-w_{il}^{pre}/I$  and submit it to the VPP operator.

Under a decentralized scheme, a consumer $i$ in SR $l$ only has prediction information about local renewable output $W_l$, which is $w_{il}^{pre}$. Its decision-making problem is
\begin{align}
\mathop{\max}_{D_{il}}~ & U_i(D_{il})-D_{il} \mathbbm{E}[\lambda(D,W)|W_{il}^{pre}=w_{il}^{pre}]
\end{align}

A candidate affine policy $D_{il}^{dis*}$ is
\begin{align}
\label{eq:individualaffineDPD}
D_{il}^{dis*}:=\hat{b}_1^l+\hat{b}_{2i}^l(w_{il}^{pre}-\overline{w}_l),\forall i \in \mathcal{I}_l, \forall l \in \mathcal{L}   
\end{align}

Applying the law of large numbers,  with many consumers in SR $l$, the average of $D_{il}^{dis*},\forall i \in \mathcal{I}_l$ is approximately
\begin{align}
\label{eq:16}
\overline D_{l}^{dis*}:=\hat{b}_1^l+\hat{b}_{2}^l(w_{l}-\overline{w}_l), \forall l \in \mathcal{L}
\end{align}
where $\hat{b}_2^l={\sum_{i \in \mathcal{I}_l} \hat{b}_{2i}^l}/{I}$. Then the optimal consumption level of consumer $i \in \mathcal{I}_l$ under given $\beta_l$ and precision $\tau_{il}$ is given by the following proposition. 
\begin{proposition}
	\label{prop3}
	Given the updates $\beta_l,\forall l \in \mathcal{L}$ and prediction precision $\tau_{il},\forall l \in \mathcal{L}, i \in \mathcal{I}_l$, the consumer $i \in \mathcal{I}_l$'s optimal demand under the decentralized scheme is given by \eqref{eq:individualaffineDPD} with
	\begin{align}
	\label{eq:dis-consumer}
	\hat{b}_{2i}^l =~& \frac{\alpha \tau_{il}}{\gamma \tau_l+u},\forall i \in \mathcal{I}_l, \nonumber\\
	\hat{b}_{1}^l =~ & \frac{t-\beta_l+\alpha \overline{w}_l}{\gamma+u}
	\end{align}
	where $\gamma=I\alpha$.
	%\vspace{1mm}
\end{proposition}

It follows that the consumer $i \in \mathcal{I}_l$'s expected net utility is
\begin{align}
\mathbbm{E}[\pi_i]:=~ & \mathbbm{E}[\frac{1}{2}u(D_{il}^{dis*})^2] \nonumber\\
=~ & \frac{1}{2}u\left[( \frac{t-\beta_l+\alpha \overline{w}_l}{\gamma+u})^2+\frac{\alpha^2\sigma_{W_l}^2\tau_{il}}{(\gamma \tau_l+u)^2}\right]
\end{align}

The proof of Proposition \ref{prop3} is similar to that of Proposition \ref{prop1} and so is omitted here. Different from  the optimal decision under centralized scheme, here the influence of renewable uncertainties in other SRs is reflected in the ``constant'' term $\hat{b}_1^l$. Meanwhile, we can figure out that the optimal prediction precision of SR $l$ is independent of the information from other SRs. This is formalized in the following proposition.
\begin{proposition}
	\label{prop4}
	With decentralized predictions, the optimal prediction precision for the SR $l$ is
	\begin{align}
	\label{eq:optimal-dis-pa}
	\tau_{l}^{dis*}=\mathop{\max}\left(0,\frac{ \sigma_{W_l}^2-\sqrt{2mu/I}/\alpha}{\sigma_{W_l}^2+ \sqrt{2mI/u}}\right),\forall l \in \mathcal{L}
	\end{align}
	\vspace{0.1mm}
\end{proposition}

The proof of Proposition~\ref{prop4} follows the same line of the proof of Proposition~\ref{prop2}. With the above analysis, we close this subsection with several remarks, indicating the benefits of our decentralized scheme:

1) The optimal value of $\tau_l^{dis*}$ is independent of $\beta_l$ and remains unchanged in each iteration. This allows each SR to determine its optimal prediction precision individually and beforehand, without communicating with the others.

2) The proposed decentralized scheme is in accord with the centralized scheme in two extreme cases: 
\begin{itemize}
	\item [$\bullet$] When $m=0$, improving prediction technique doesn't incur any cost. At this time, $\tau_l^{cen*}=\tau_l^{dis*}=1$, which means the more accurate the better.
	\item [$\bullet$] When $m$ is sufficiently large, improving prediction precision is so costly such that consumers choose not to improve prediction techniques, i.e., $\tau_l^{cen*}=0$ or $\tau_l^{dis*}=0$.
\end{itemize}
\subsection{Convergence}
It has been shown in Proposition \ref{prop4} that the optimal prediction precision can be determined individually by each SR. However, the optimal demand level still depends on the information from other SRs. In this subsection, a condition under which the DPP algorithm (\textbf{Algorithm 1}) converges is given. The optimal decision at convergence point is also revealed. Consider the following inequality:

\begin{enumerate}[label=C\arabic*:]
	\item  $(L-2)I\alpha < u$.
\end{enumerate}

\begin{proposition}
	\label{prop5}
	When C1 holds and when every consumer in every SR $l$ chooses the optimal prediction precision $\tau_{l}^{dis*}$ as in \eqref{eq:optimal-dis-pa}, the average demand of every SR $l\in\mathcal{L}$ in  the DPP algorithm converges to
	\begin{align}
	\label{eq:opt-d-dis}
	\hat{D}_l^{dis*} =~ & \frac{t-\beta_0+\alpha \sum \limits_{j=1}^L \overline{w}_j}{\zeta+u} + \frac{\alpha (\tau_{l}^{dis*}-1)}{\gamma \tau_{l}^{dis*}+u} (w_l-\overline{w}_l) \nonumber\\
	~& + \frac{\gamma+u}{\zeta + u} \sum \limits_{j=1}^L \frac{\alpha}{\gamma \tau_{j}^{dis*}+u}( w_j - \overline{w}_j)
	\end{align}
\end{proposition}

The proof of Proposition \ref{prop5} can be found in Appendix~\ref{apen-prop5}. It indicates that even if the consumer in SR $l$ has no access to the predictions of the renewable outputs from other SRs, its optimal decision, however, converges to a quantity that contains such information, if Condition C1 holds.
\begin{figure}[h]
	\centering
	\includegraphics[width=0.8\columnwidth]{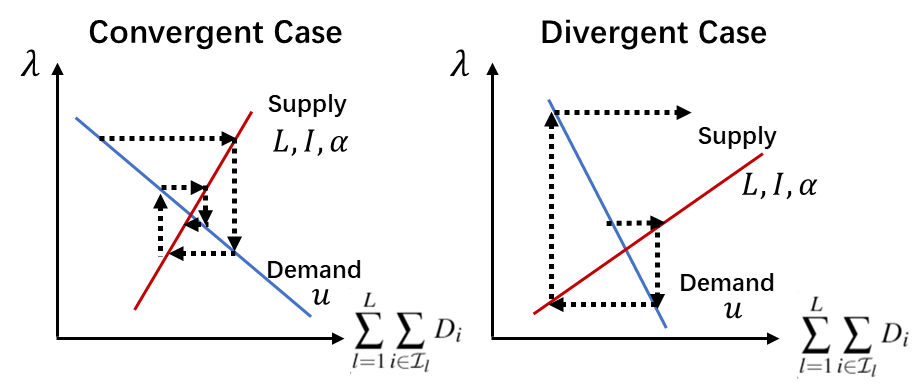}
	\caption{Economic intuition behind Condition C1.}
	\label{fig:conditionC1}
\end{figure}
%	\vspace{-2mm}
	
\textbf{Remark:} We try to give an economic interpretation of Condition C1: Recall that equation \eqref{eq:marketprice} shows the relationship of electricity price $\lambda$ and total demand $\sum_{l=1}^L \sum_{i \in \mathcal{I}_l} D_i$. A higher $\lambda$ is associated with more demand, which in other words means the market is more willing to supply. Equation \eqref{eq:marketprice} is actually the supply function of the electricity market and its price-to-quantity sensitivity for each consumer is related to $L$, $I$ and $\alpha$. In equation \eqref{eq:condition}, the demand $D_i^*$ decreases with the price $\mathbbm{E}[\lambda(D,W)|\mathcal{A}_i]$, which can be interpreted as the demand function of an individual consumer with the absolute value of the slope equals to $u$. Condition C1 means that the absolute value of the slope of the supply curve is less than the absolute value of the slope of the demand curve. When C1 holds, the deviations of price and quantity are compressed after the market reactions. A illustrative diagram is given in Fig.\ref{fig:conditionC1} and similar phenomenon can be found in \cite{zhao2014convergence}, which is also known as cobweb model in economics.

%\textcolor{blue}{Equation \eqref{eq:marketprice} is the supply function of the distribution electricity market} and its price-to-quantity sensitivity for each consumer is related to $L$, $I$ and $\alpha$. Equation \eqref{eq:condition} is the demand function of an individual consumer with the absolute value of the slope equals to $u$. The economic interpretation of Condition C1 is that the absolute value of the slope of the supply curve is less than the absolute value of the slope of the demand curve. \textcolor{blue}{When C1 holds, the deviations of price and quantity are compressed after the market reactions.} Similar phenomenon can be found in \cite{zhao2014convergence}, which is also known as cobweb model in economics.

%\textcolor{blue}{Equation \eqref{eq:marketprice} shows the relationship of electricity price $\lambda$ and total demand $\sum_{l=1}^L \sum_{i \in \mathcal{I}_l} D_i$. The higher $\lambda$, the more demand the market is willing to supply, which in other words can be interpreted as the supply function of the electricity market. In equation \eqref{eq:condition}, the demand $D_i^*$ decreases with the price $\mathbb{E}[\lambda(D,W)|\mathcal{A}_i]$, which can be interpreted as the demand function of an individual consumer.}

\subsection{Comparison: centralized v.s. decentralized}
Given the optimal decisions of consumers with both centralized and decentralized predictions, in this section, we compare the gap between them to theoretically validate the effectiveness of the DPP algorithm.

When Assumptions A1-A2 and Condition C1 hold, we are able to compute the optimal average demand for each SR, with both centralized and decentralized predictions. Denote them by $D_l^{cen*}$ (given by \eqref{eq:opt-d-cen}) and $\hat{D}_l^{dis*}$ (given by \eqref{eq:opt-d-dis}), respectively.
\begin{theorem}
	\label{prop6}
	For any $l\in\mathcal{L}$, the expectation
	\begin{align}
	\mathbbm{E}[D_l^{cen*}-\hat{D}_l^{dis*}]=0 
	\end{align}
	and the variance $\mathbbm{D}(D_l^{cen*}-\hat{D}_l^{dis*})$ is bounded.
	%\vspace{0.5mm}
\end{theorem}

The proof of Theorem~\ref{prop6} is postponed to Appendix~\ref{apen-prop6}. Theorem \ref{prop6} implies that for each SR $l$, the DPP algorithm can achieve the same expected average demand level as when centralized predictions are available. The variance of the gap between $D_l^{cen*}$ and $\hat D_l^{dis*}$ is bounded, implying that the deviation from the optimal demand with centralized predictions  incurred by information loss would likely be well controlled. With above, by providing a more practical local prediction paradigm, the DPP algorithm allows each SR to operate in a decentralized manner, while achieving an acceptable overall market performance.

\section{Simulation}\label{sec:simulation}

\subsection{Illustrative example}
In this section, numerical experiments are conducted to support the theoretical results and provide some insights. First, a simple case with 3 subregions is tested, i.e., $L=3$. The parameters for each SR are: $\overline{w}_1=300$ kWh, $\overline{w}_2=500$ kWh, $\overline{w}_3=400$ kWh, and $\sigma_{W_1}^2=400 \;\mbox{(kWh)}^2$, $\sigma_{W_2}^2=3000\; \mbox{(kWh)}^2$, $\sigma_{W_3}^2=1600 \;\mbox{(kWh)}^2$. We choose $m=5 \; \mbox{\$(kWh)}^2$, $u=8 \; \mbox{\$/(kWh)}^2$, $t=80 \; \mbox{\$/(kWh)}$, $\alpha=0.003  \; \mbox{\$/(kWh)}^2$ and $\beta_0=0.03 \; \mbox{\$/(kWh)}$. 
%\begin{table}[h]
%	\renewcommand{\arraystretch}{1.3}
%	\centering
%	\caption{Parameters of %subregions}
%	\label{tab:para}
%	\begin{tabular}{ccccc}
%		\hline 
%		$l$ & 1 & 2 & 3 \\
%		\hline
%		$\overline{w}_l / \mbox{kWh}$ %& 300 & 500 & 400 \\
%		$\sigma_{W_l}^2 / %\mbox{(kWh)}^2$ & 400 & 3000 & 1600 %\\ 
%		\hline
%	\end{tabular}
%\end{table}

\begin{figure}[h]
	\centering
	\includegraphics[width=1.0\columnwidth]{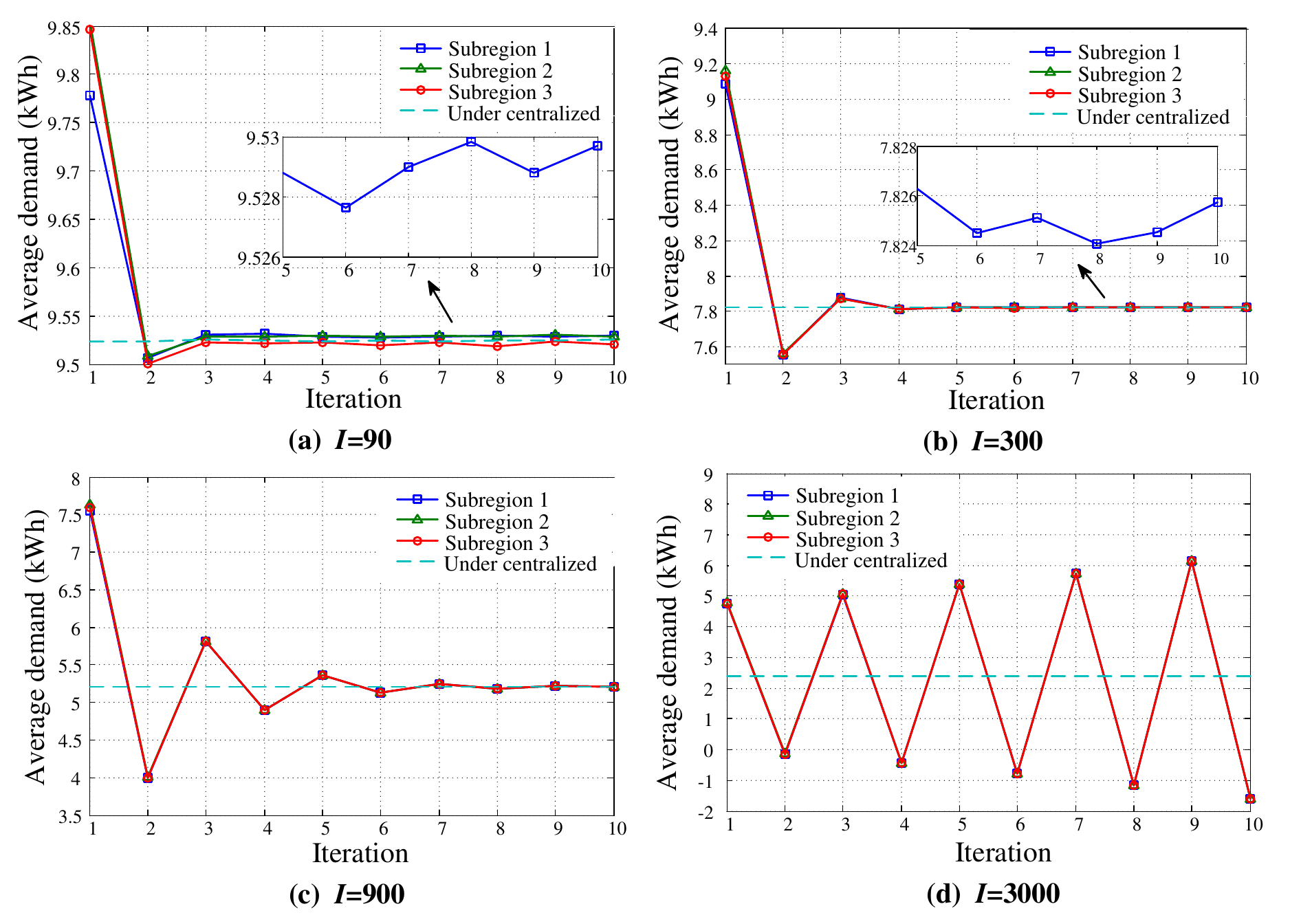}
	%	\vspace{-2mm}
	\caption{Performance of the DPP algorithm. When Condition C1 holds, the average demand converges; otherwise it does not.}
	\label{fig:convergence}
\end{figure}

First, the performance of the DPP algorithm with different numbers of consumers is investigated. We test the cases with $I=90$, $I=300$, $I=900$, and $I=3000$, and the results are shown in Fig.~\ref{fig:convergence}.(a)-(d), respectively. Note that Condition C1 holds for cases (a)-(c), and thus, the DPP algorithm converges as shown in the figure. However, the DPP algorithm fails to converge for case (d) because $(L-2)I\alpha = 9>8=u$. This verifies that the proposed algorithm is guaranteed to converge with a moderate number of consumers. As proved in Proposition \ref{prop5}, the DPP algorithm converges when Condition C1 holds. In other words, given the number of subregions $L$, parameters $\alpha$ and $u$, the number of consumers needs to satisfy $I<u/\left(\alpha(L-2)\right)$. Under above settings the DPP algorithm converges when $I<2666$, which is sufficient for most cases of practical interests. When there is an extremely large number of consumers, we can aggregate some of them as decision-making coalitions.

Additionally, we observe that with more consumers, the gap between centralized and DPP algorithms at the equilibrium point is narrower, while the convergence becomes slower with severer oscillation (which can also be verified in the proof of Proposition \ref{prop5}, by noticing that the spectral radius of $\textbf{H}$ is monotonically increasing with consumer number $I$), revealing a possible trade-off between transient and steady-state performance of DPP.
%more participants can benefit the algorithm to some extent. In short,
%	\vspace{-2mm}

\begin{figure}[h]
	\centering
	\includegraphics[width=1.0\columnwidth]{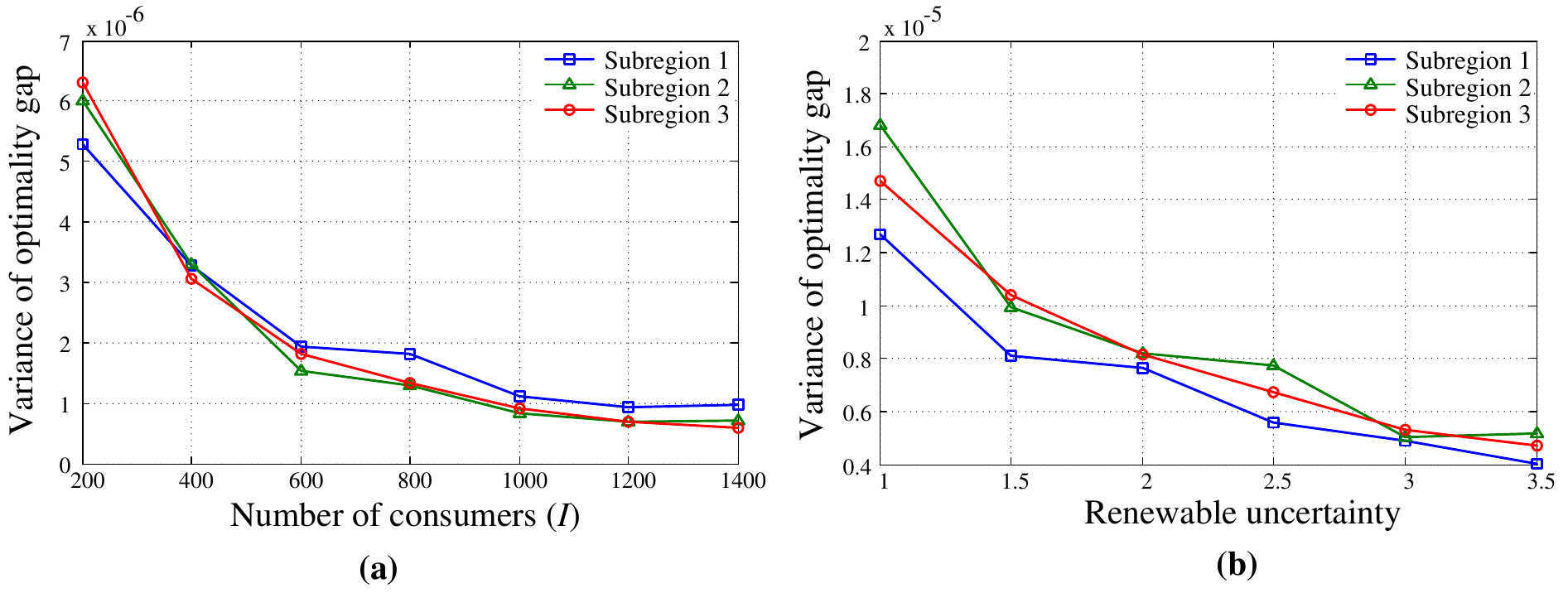}
	\caption{Impact of the number of consumers and uncertainty on the variance.}
	\label{fig:influence}
\end{figure}

We further analyze the impact of consumer number $I$ and uncertainty level $\sigma_{W_l}^2$ on the gap between $D_l^{cen*}$ and $\hat D_l^{dis*}$. We change $I$ from $200$ to $1400$, and $\sigma_{W_l}^2$ from $1$ to $3.5$ times of their original values, for all $l \in \mathcal{L}$. With these parameters changing, the variance $\mathbbm{D}(D_l^{cen*}-\hat{D}_l^{dis*})$ is plotted in Fig.~\ref{fig:influence}.(a)-(b), respectively. The gaps decline with the growth of either factor, implying that the DPP algorithm can perform better with more consumers and higher uncertainty. Even with a small number of consumers $I$ and low uncertainty $\sigma_{W_l}^2,\forall l \in \mathcal{L}$, the gaps are still small. In fact, two factors account for the gap between $D_l^{cen*}$ and $\hat D_l^{dis*}$: one is the gap between optimal prediction accuracy $\tau_l^{cen*}$ given in \eqref{eq:prediction-accuracy-cen} and $\tau_l^{dis*}$ given in \eqref{eq:optimal-dis-pa}. This gap decreases with growing $I$ or $\sigma_{W_l}^2$. The other factor is the approximation of $w$ in Algorithm 1 using the average of predictions $w_{il}$. In the centralized scheme, $I L$ predictions are summed up to approximate $w$, while in the decentralized scheme, only $I$ predictions are used. As $I$ increases, the two schemes tend to the same approximation accuracy. Considering the trends of these two factors, the variance $\mathbbm{D}(D_l^{cen*}-\hat D_l^{dis*})$ declines with the growth of either $I$  or renewable uncertainty, as shown in Fig.\ref{fig:influence}.

\begin{figure}[h]
	\centering
	\includegraphics[width=0.95\columnwidth]{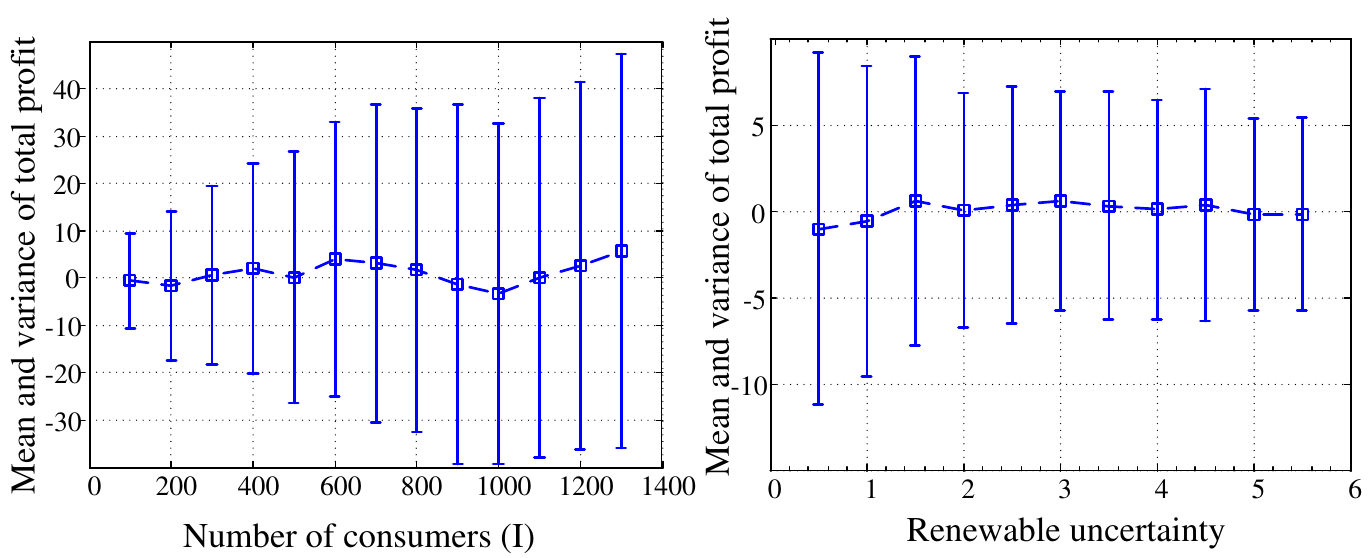}
	\caption{Influence of the number of consumers and uncertainty on the  profit.}
	\label{fig:influence-profit}
\end{figure}

Besides the optimal consumption strategy above, we also analyze the expected total profit, defined as the following:.
\begin{align}
 \nonumber
 ~ &   \mathbbm{E}\left[\sum_{l=1}^L \sum_{i \in \mathcal{I}_l} U_i(D_{i})-\lambda(D,W)D\right]-\sum_{l=1}^L \frac{m}{\sigma_{W_l}^2}\frac{\tau_{l}}{1-\tau_{l}}
\end{align}
We calculate the total profit of the centralized method minus that of the decentralized one, and show the mean and variance of this difference in Fig. \ref{fig:influence-profit}, with a varying number of  consumers $I$ and uncertainty level $\sigma_{W_l}^2$.

We observe that the expected total profit of the centralized method does not necessarily exceed that of the decentralized one. We briefly discuss our conjecture about the reason, while a formal analysis is left for future work. 
Under the proposed settings, two main factors that may influence the total profit are interest conflicts among consumers (which have been widely investigated in the literature based on the Cournot model) and renewable uncertainty. 
When consumers do not know the exact value of market price and hence decide their strategies by predictions, the associated ambiguity may offset the negative impact of interest conflicts. If this happens, the expected total profit of the decentralized method may surpass that of the centralized one. 
As $I$ increases, the impact of interest conflicts starts to dominate, making the centralized method steadily outperform the decentralzied one. 
From another perspective, as $\sigma_{W_l}^2$ increases, the gap between centralized and decentralized methods wanders near zero since uncertainty dominates and diminishes the effect of interest conflicts.

\begin{table}[h]
	\renewcommand{\arraystretch}{1.3}
	\centering
	\caption{Time and number of iterations under different $L$}
	\label{tab:time}
	\begin{tabular}{ccccccc}
		\hline 
		$L$ & 3 & 9 & 15 & 21 & 27 & 33 \\
		\hline
		Time (s) & 0.078 & 0.157 & 0.231 & 0.259 & 0.603 & 0.428 \\
		Iterations & 3 & 6 & 9 & 10 & 17 & 71 \\ 
		\hline
	\end{tabular}
\end{table}
To show the practicability of the proposed algorithm, we test the cases with more subregions by letting $L$ equal to 3, 9, 15, 21, 27, and 33, respectively. The computational time and number of iterations needed to converge are shown in TABLE \ref{tab:time}. Though both time and iterations grow with the scale of the test system, they are all moderate and acceptable.

\begin{figure}[h]
	\centering
	\includegraphics[width=0.75\columnwidth]{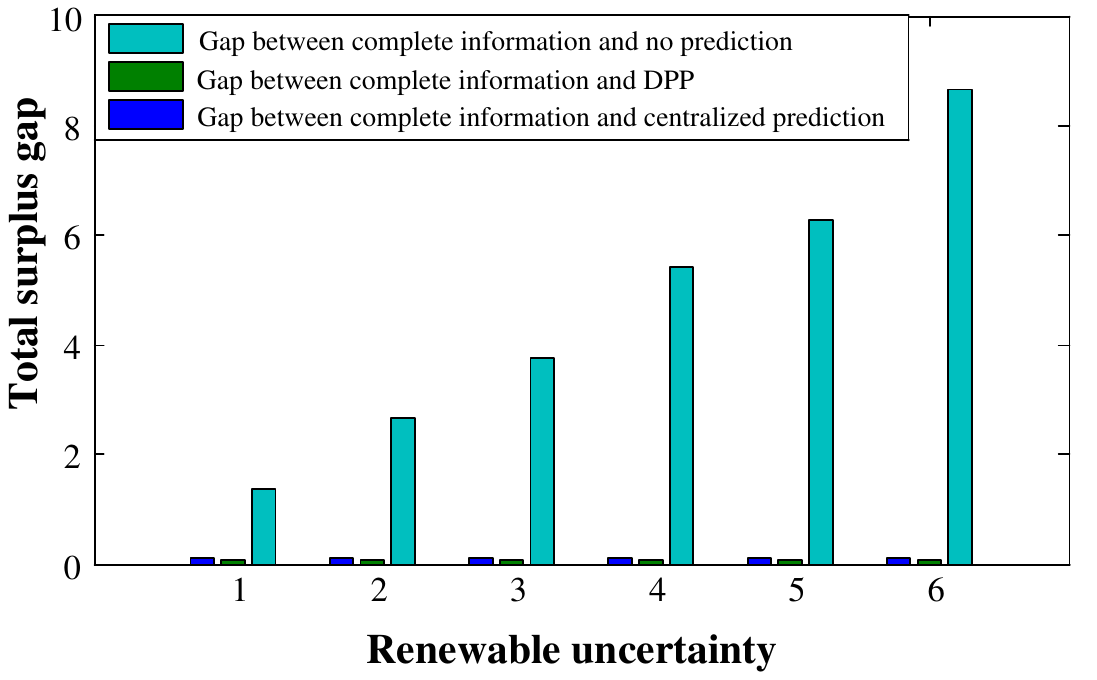}
	\caption{Total surplus gaps between complete information and no prediction, DPP, centralized prediction.}
	\label{fig:total-surplus}
\end{figure}

To further show the impact of proposed algorithm on market efficiency, we change the renewable uncertainty variance from 1 to 6 times of its original value, and record the total surplus gaps as in Fig.\ref{fig:total-surplus}. We can observe that the gaps between complete information \& centralized prediction and between complete information \& DPP are almost zero, showing that the proposed method can achieve near-optimal market efficiency.

\subsection{Practical example}
We further extend the model \eqref{eq:consumer} to incorporate multiple periods, capacity constraint, and total load requirement:
\bsq
\begin{align}
    \mathop{\max}_{D_{it},\forall t=1,...,T} ~ & \pi_i:=\sum \limits_{t=1}^T \left(U_i(D_{it})-D_{it}\mathbbm{E}[\lambda_t(D_t,W_t)|\mathcal{A}_i]\right) \\
    \mbox{s.t.}~ & D_{i}^{down} \le D_{it} \le D_{i}^{up}, \forall t=1,...,T \\
    ~ & \sum \limits_{t=1}^T D_{it}=D_{i}^{total}
\end{align}
\esq
where $\lambda_t(D_t,W_t)$ at each period   $t$ is given by \eqref{eq:marketprice} with $D_{it},\forall i \in \mathcal{I}$. $D_{i}^{down}$ and $D_{i}^{up}$ are the lower and upper bound of the capacity limit, respectively. $D_{i}^{total}$ is the total demand requirement of consumer $i \in \mathcal{I}$ over $T$ periods. With complete information, each consumer knows exactly the renewable output $w_t$, so the objective function becomes
\begin{align}
    \mathop{\max}_{D_{it},\forall t=1,...,T}~ \pi_i:=\sum \limits_{t=1}^T \left(U_i(D_{it})-D_{it}\lambda_t(D_t,w_t)\right)
\end{align}

With no prediction, each consumer makes decision according to the expected value $\bar w_t$, and the objective function becomes
\begin{align}
    \mathop{\max}_{D_{it},\forall t=1,...,T}~ \pi_i:=\sum \limits_{t=1}^T \left(U_i(D_{it})-D_{it}\lambda_t(D_t,\bar w_t)\right)
\end{align}

\begin{figure}[h]
	\centering
	\includegraphics[width=0.75\columnwidth]{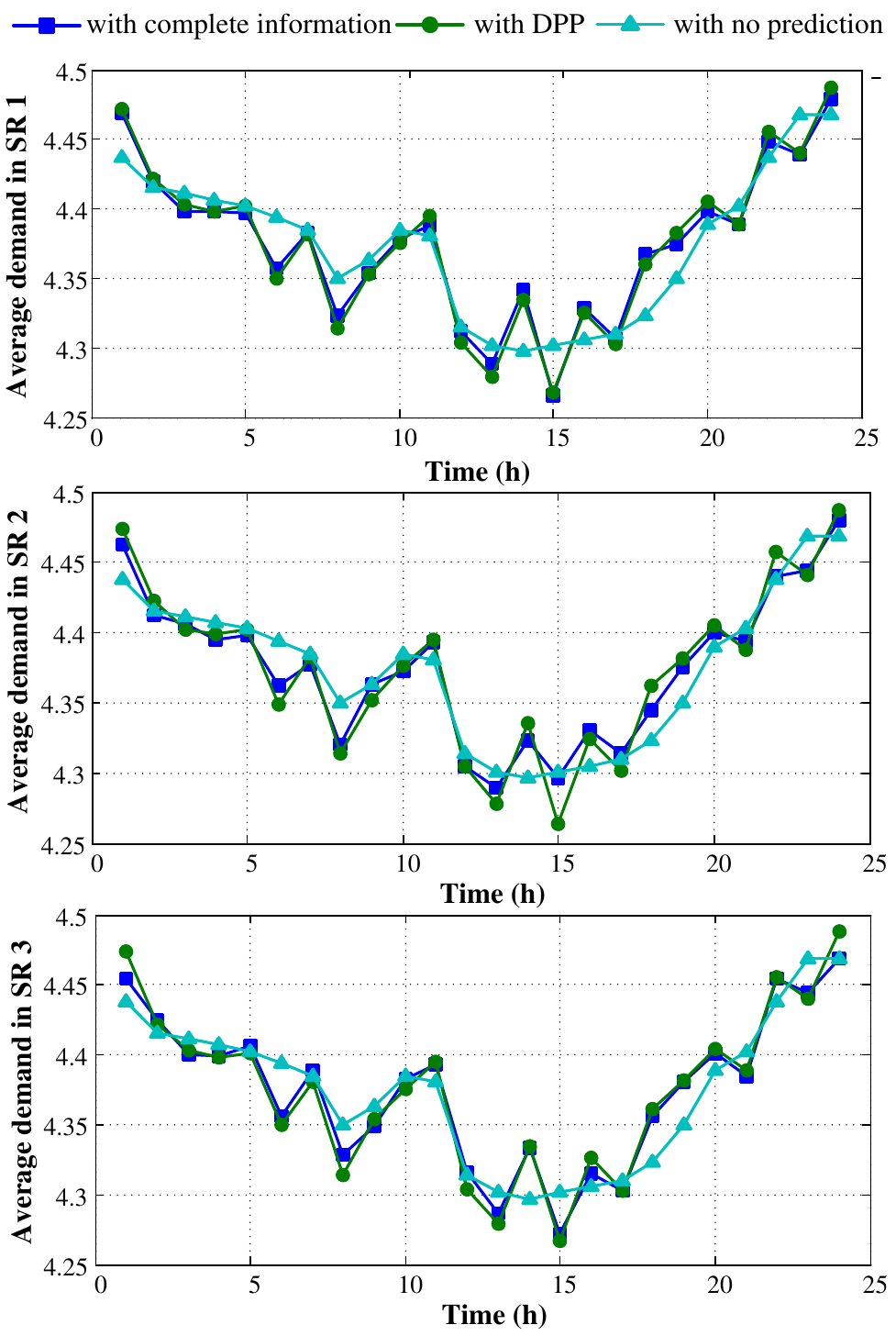}
	\caption{Demand under complete information, DPP and no prediction..}
	\label{fig:practical-case1}
\end{figure}

\begin{figure}[h]
	\centering
	\includegraphics[width=0.75\columnwidth]{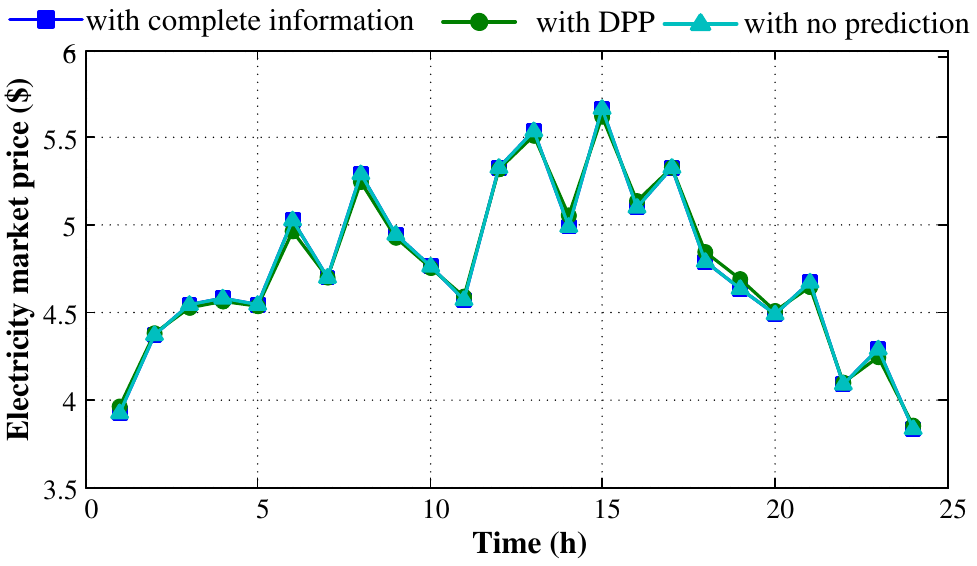}
	\caption{Electricity prices under complete information, DPP and no prediction.}
	\label{fig:practical-case2}
\end{figure}

We simulate consumers' optimal demand with complete information, DPP, and no prediction, respectively. Let $D_i^{total}=105$ kWh with $D_{it}^{down}=2$ kWh and $D_{it}^{up}=5$ kWh. The prediction accuracy under DPP is set to $\tau_l^{dis}=0.5$. Other parameters are the same as the benchmark case. The consumers’ demands are recorded in Fig.\ref{fig:practical-case1} and the real-time electricity prices are displayed in Fig.\ref{fig:practical-case2}. We can find that with DPP, consumers' optimal demand profiles are closer to the curves under complete information, while the curves without prediction differ a lot. Moreover, consumers’ demand with DPP is complementary to the electricity market prices, i.e. when the price is higher a consumer chooses to use less electricity, and vice versa. This verifies that our proposed algorithm can provide sufficient information for enhancing consumer’s profit.

\section{Conclusion}
Virtual power plant (VPP) can gather consumers from different geographical locations and help improve their market profits. Facing the increasing penetration of distributed renewables, one desired service of the VPP is to provide predictions of uncertain renewable outputs to its consumers. These predictions can be provided either in a centralized or a decentralized manner. Since the former one may encounter the problems of data privacy and computational burden, a decentralized prediction provision (DPP) algorithm is proposed in this paper, enabling consumers to decide their consumption levels with local predictions. We have proved that providing predictions can boost the social total surplus. Several appealing properties of the proposed algorithm, including its convergence guarantee and statistical performance compared to the centralized framework, have been proved. Illustrative examples validate the theoretical results and provide some further insights: the variance of the optimal demand gap between centralized and decentralized schemes decreases with more consumers and higher level of uncertainty; the total profit of consumers exhibits the trade-off between interest conflicts and renewable uncertainties. A practical case is also tested to show scalability of the proposed method. Future research entails quantifying the gap of total profit under centralized and decentralized schemes, as well as generalizing the proposed algorithm to VPPs with diverse facilities.

	\ifCLASSOPTIONcaptionsoff
	\newpage
	\fi
	
	\bibliographystyle{IEEEtran}
	\bibliography{IEEEabrv,mybib}

\appendix
\makeatletter
\@addtoreset{equation}{section}
\@addtoreset{theorem}{section}
\makeatother
\setcounter{equation}{0}  
\renewcommand{\theequation}{A.\arabic{equation}}
\subsection{Proof of Lemma \ref{lemma1}}
\label{apen-lemma1}
Note that $\{\varepsilon_{il}:i\in\mathcal{N}_{l}\}$ are independent and identically distributed (i.i.d.) and $\varepsilon_{il}$ is independent from $W_l$.
This implies that $\mathbbm{E}[\varepsilon_{il}|W_l=w_l]=\mathbbm{E}[\varepsilon_{il}]=0$ and $\mathrm{\mathrm{cov}}(\varepsilon_{il},W_l)=0$ for each subregion $l\in\mathcal{L}$ and consumer $i\in\mathcal{N}_l$. Therefore, $\mathbbm{E}[W_{il}^{pre}|W_l=w_l]=w_l$. 
We have:
%{\color{ForestGreen}In the following equations, I changed $w_{il}^{pre}$ to $W_{il}^{pre}$ whenever it is a random variable instead of a specific realization.}
\begin{align}
\overline{w}_l = \mathbbm{E}[W_l] = \mathbbm{E}[\mathbbm{E}[W_l|W_{il}^{pre}]]= A_{1i,l} +A_{2i,l} \overline{w}_l
\end{align}
and with $\mathrm{cov}(\varepsilon_{il},W_l)=0$, it gives
\begin{align}
\sigma_{W_l}^2=~ & \mathrm{cov}(W_l,W_{il}^{pre}) \nonumber\\
= ~ & \mathrm{cov}(W_{il}^{pre},\mathbbm{E}[W_l|W_{il}^{pre}]) \nonumber \\
= ~ &\mathrm{cov}(W_{il}^{pre},A_{1i,l}+A_{2i,l}W_{il}^{pre}) \nonumber\\
=~ & A_{2i,l}(\sigma_{W_l}^2+\sigma_{\varepsilon_{il}}^2).
\end{align}

With the above two conditions, we obtain
\bq
A_{2i,l} = \frac{\sigma_{W_l}^2}{\sigma_{W_l}^2+\sigma_{\varepsilon_{il}}^2}, \;\; A_{1i,l}=(1-A_{2i,l})\overline{w}_l.
\eq

\setcounter{equation}{0}  
\renewcommand{\theequation}{B.\arabic{equation}}
\subsection{Proof of Proposition \ref{prop1}}
\label{apen-prop1}
Substitute \eqref{eq:consumer-stratetgy} into the optimal condition \eqref{eq:condition}, we get
\begin{align}
~ & t-u(b_1+\sum \limits_{l=1}^L b_{2i}^l(w_{il}^{pre}-\overline{w}_l)) \nonumber\\
=~& \mathbbm{E}[\zeta (b_1+\sum \limits_{l=1}^L b_2^l(W_l-\overline{w}_l)) 
-\alpha \sum \limits_{l=1}^L W_l + \beta_0|  \mathcal{A}_i] \nonumber\\
=~& \mathbbm{E}[\zeta (b_1-\sum \limits_{l=1}^L b_2^l\bar{w}_l) +\beta_0 +\sum \limits_{l=1}^L (\zeta b_2^l-\alpha)W_l| \mathcal{A}_i] \nonumber\\
=~& \zeta (b_1-\sum \limits_{l=1}^L b_2^l\overline{w}_l) + \beta_0 + \sum \limits_{l=1}^L (\zeta b_2^l - \alpha ) (	1-\tau_{il})\overline{w}_l \nonumber\\
~ & + \sum \limits_{l=1}^L (\zeta b_2^l-\alpha) \tau_{il} w_{il}^{pre}.
\end{align}

This must hold for all $w_{il}^{pre}$, whence
\begin{align}
t-u(b_1-\sum \limits_{l=1}^L b_{2i}^l \overline{w}_l) =~ & \zeta (b_1-\sum \limits_{l=1}^L b_2^l\overline{w}_l) + \beta_0 \nonumber\\
~ & + \sum \limits_{l=1}^L (\zeta b_2^l - \alpha ) (1-\tau_{il}) \overline{w}_l \nonumber,\\
\text{and } -u b_{2i}^l  =~ & (\zeta b_2^l - \alpha )\tau_{il},\forall l \in \mathcal{L}.
\end{align}

With above equations, we can get the values of $b_1$ and $b_{2i}^l,\forall l \in \mathcal{L}, i \in \mathcal{N}_l$:
\begin{align}
b_{2i}^l = \frac{\alpha \tau_{il}}{\zeta \tau_{l} + u},\forall l \in \mathcal{L}, i \in \mathcal{N}_l\; , \; b_1 =  \frac{t-\beta_0+\sum_l \alpha \overline{w}_l}{\zeta+u}.
\end{align}

\setcounter{equation}{0}  
\renewcommand{\theequation}{C.\arabic{equation}}
\subsection{Proof of Proposition \ref{prop2}}
\label{apen-prop2}
For the centralized prediction scheme, the prediction cost for each consumer $i \in \mathcal{I}$ is $\sum_{l \in \mathcal{L}} h_l(\tau_{il})$ with $\check{m}=m/(IL)$. Then
\begin{align}
\frac{\partial \left(\mathbbm{E}[\pi_{i}]-\sum \limits_{l=1}^L h_{l}\right)}{\partial \tau_{il}}=  \frac{u\alpha^2\sigma_{W_l}^2}{2(\zeta \tau_{l}+u)^2}-\frac{m}{IL\sigma_{W_l}^2}\frac{1}{(1-\tau_{il})^2}.
\end{align}

Because of symmetry, $\tau_{il}=\tau_l,\forall i \in \mathcal{N}_l=\mathcal{I}$. Therefore, the optimal prediction precision is
\begin{align}
\tau_{l}^{cen*}=\mathop{\max}\left(0,\frac{ \sigma_{W_l}^2-\sqrt{2mu/(IL)}/\alpha}{\sigma_{W_l}^2+ \sqrt{2mIL/u}}\right),\forall l \in \mathcal{L}.
\end{align}

\setcounter{equation}{0}  
\renewcommand{\theequation}{D.\arabic{equation}}
\subsection{Proof of Theorem \ref{prop7}}
\label{apen-prop7}
First, we calculate the value $\mathbbm{E}[TS(D^0,W)]$. Similar to \eqref{eq:condition}, we have $D^0$ satisfies
\begin{align}
\label{eq:first-best-condition}
    t-uD^0 = \alpha ILD^0-\alpha \sum \limits_{l=1}^L w_l + \beta_0
\end{align}
With assumption A1, definition \eqref{eq:total-surplus} is equivalent to
\begin{align}
    TS(D,W)=ILU(D)-\int_0^{ILD} \lambda(x,W)dx
\end{align}
Then
\begin{align}
    ~ & TS(D^0,W) \nonumber\\
    =~ &  IL\left(-\frac{u}{2}(D^0)^2+tD^0 - \frac{\alpha}{2}IL(D^0)^2+\alpha D^0 \sum \limits_{l=1}^L w_l - \beta_0 D^0\right) \nonumber\\
    = ~ & IL \frac{u+\zeta}{2}(D^0)^2
\end{align}
From \eqref{eq:first-best-condition}, we have
\begin{align}
    D^0=\frac{t+\alpha \sum \limits_{l=1}^L w_l -\beta_0}{u+\zeta} = b_1 + \frac{\alpha \sum \limits_{l=1}^L (w_l-\overline{w}_l)}{u+\zeta}
\end{align}
Therefore
\begin{align}
    \mathbbm{E}[TS(D^0,W)]=IL\frac{u+\zeta}{2}\left[b_1^2+\left(\frac{\alpha}{u+\zeta}\right)^2\sum \limits_{l=1}^L \sigma_{W_l}^2\right]
\end{align}

Next, we calculate $\mathbbm{E}[\mathbbm{E}[TS(D^{cen*},W)|W=w]]$.
\begin{align}
\label{eq:with-predict-expected}
    ~ & \mathbbm{E}[TS(D^{cen*},W)|W=w] \nonumber\\
    =~ & \sum \limits_{l=1}^L \sum \limits_{i \in \mathcal{I}_l} U(D_i^{cen*})-\int_{0}^{IL\overline{D}^{cen*}}\lambda(x,W)dx
\end{align}
where
\begin{align}
    D_i^{cen*}=~ & b_1+\sum \limits_{l=1}^L b_{2i}^l (w_{il}^{pre}-\overline{w}_l)\nonumber\\
    \overline{D}^{cen*}=~ & b_1+\sum \limits_{l=1}^L b_2^l(w_l-\overline{w}_l) \nonumber
\end{align}
Then the expectation of the first term in \eqref{eq:with-predict-expected} is
\begin{align}
    ~ & \mathbbm{E}[\sum \limits_{l=1}^L \sum \limits_{i \in \mathcal{I}_l} U(D_i^{cen*})] \nonumber\\
    =~ & -\frac{uIL}{2} \left(b_1^2+\sum \limits_{l=1}^L (b_2^l)^2 \left(\sigma_{W_l}^2+\sigma_{\epsilon_{l}}^2\right)\right)+ILtb_1
\end{align}
The expectation of the second term in \eqref{eq:with-predict-expected} is
\begin{align}
    ~ & \mathbbm{E}[\int_{0}^{IL\overline{D}^{cen*}}\lambda(x,W)dx] \nonumber\\
    =~ & \frac{\zeta IL}{2} \left(b_1^2+\sum \limits_{l=1}^L (b_2^l)^2 \sigma_{W_l}^2\right)+\beta_0 IL b_1  \nonumber\\
    ~ & -\zeta (\sum \limits_{l=1}^L \overline{w}_l)b_1-\zeta \sum \limits_{l=1}^L b_2^l \sigma_{W_l}^2
\end{align}
Therefore
\begin{align}
\label{eq:total-surplus-cen}
    ~& \mathbbm{E}[\mathbbm{E}[TS(D^{cen*},W)|W=w]] \nonumber\\
    =~& IL\frac{u+\zeta}{2}b_1^2 +\frac{\alpha IL }{2}\sum \limits_{l=1}^L b_2^l \sigma_{W_l}^2 
\end{align}
The first term in \eqref{eq:total-surplus-cen} is constant. The second term increases with $b_2^l,\forall l \in \mathcal{L}$ and thus increases with $\tau_{l}$ until $\tau_l=1,\forall l \in \mathcal{L}$, at which the case with improved predictions degenerates to the case with complete information:
\begin{align}
\mathbbm{E}[TS(D^0,W)]=\mathbbm{E}[\mathbbm{E}[TS(D^{cen*},W)|W=w]]|_{\tau_l=1,\forall l\in \mathcal{L}} \nonumber    
\end{align}
When $\tau_{l}$ decreases to $\tau_l=0,\forall l \in \mathcal{L}$,  it degenerates to the case without improved predictions:
\begin{align}
\mathbbm{E}[TS(D^1,W)]=\mathbbm{E}[\mathbbm{E}[TS(D^{cen*},W)|W=w]]|_{\tau_l=0,\forall l\in \mathcal{L}} \nonumber    
\end{align}
We can draw the conclusion that
\begin{align}
    \mathbbm{E}[TS(D^0,W)] \ge \mathbbm{E}[\mathbbm{E}[TS(D^{cen*},W)|W=w]] \ge \mathbbm{E}[TS(D^1,W)] \nonumber
\end{align}

\setcounter{equation}{0}  
\renewcommand{\theequation}{E.\arabic{equation}}
\subsection{Proof of Proposition \ref{prop5}}
\label{apen-prop5}
Denote by $D_{l}^{k+1}$ the average of $D_{il}^{k+1}, \forall i \in \mathcal{I}_l$. Recall that
\begin{align}
F_l^{k+1}\approx ~ & I D_l^{k+1}-w_l \nonumber\\
=~ &  \frac{t-\beta_l^k+\alpha \overline{w}_l}{\alpha+u/I} + \frac{\alpha \tau_{l}^{dis*}}{\alpha \tau_l^{dis*}+u/I} (w_l-\overline{w}_l)-w_l \nonumber
\end{align}
and 
$
\beta_l^k = \beta_0 + \alpha \sum_{j \ne l} F_j^{k},\forall l \in \mathcal{L}.$
Thus, we have
\begin{align}
\begin{bmatrix}
F_1^{k+1} \\  \vdots \\F_L^{k+1}
\end{bmatrix}
=
\textbf{H}
\begin{bmatrix}
F_1^{k} \\  \vdots \\F_L^{k} 
\end{bmatrix}
+ \textbf{h}
\end{align} 
where
\begin{align}
\textbf{H} = ~ & \begin{bmatrix}
0 & -\frac{\gamma}{\gamma+u} & \cdots & -\frac{\gamma}{\gamma+u} \\
-\frac{\gamma}{\gamma+u} & 0 & \cdots & -\frac{\gamma}{\gamma+u} \\
\vdots & \vdots & \ddots & \vdots \\
-\frac{\gamma}{\gamma+u} & -\frac{\gamma}{\gamma+u} & \cdots & 0
\end{bmatrix} \nonumber\\
\textbf{h} = ~ & \begin{bmatrix}
\frac{t-\beta_0+\alpha \overline{w}_1}{\alpha+u/I}+\frac{\alpha \tau_{1}^{dis*}}{\alpha \tau_1^{dis*}+u/I} (w_1-\overline{w}_1)-w_1 \\
\vdots \\
\frac{t-\beta_0+\alpha \overline{w}_L}{\alpha+u/I}+\frac{\alpha \tau_{L}^{dis*}}{\alpha \tau_L^{dis*}+u/I} (w_L-\overline{w}_L)-w_L
\end{bmatrix}.
\end{align}

Since the eigenvalues of $\textbf{H}$ are $\frac{\gamma}{\gamma+u}$ and $\frac{\gamma(1-L)}{\gamma+u}$. When C1 holds, the spectral radius of $\textbf{H}$ is less than 1, so the DPP algorithm converges. Let $F_l^{k+1} = F_l^{k}=F_l^*$ and $\hat D_l^{dis*} = (F_l^*+w_l)/I$. It is easy to obtain the value of $\hat D_l^{dis*}$.
\setcounter{equation}{0}  
\renewcommand{\theequation}{E.\arabic{equation}}
\subsection{Proof of Theorem \ref{prop6}}
\label{apen-prop6}
The expectation $\mathbbm{E}[D_l^{cen*}-\hat{D}_l^{dis*}]=0$ by definition. It remains to prove that the variance is bounded. First, $\mathbbm{E}[D_l^{cen*}-\hat{D}_l^{dis*}]=0$  implies  $\mathbbm{D}(D_l^{cen*}-\hat{D}_l^{dis*})=\mathbbm{E}[(\hat{D}_l^{dis*}-D_l^{cen*})^2]$. Continuing from this, the expectation $\mathbbm{E}[(\hat{D}_l^{dis*}-D_l^{cen*})^2]$ is
%\vspace{0.5mm}
\begin{align}
\label{eq:variance-opt}
~&  \sum_{j \in \mathcal{L}, j \ne l}\underbrace{\left(\frac{\gamma+u}{\zeta + u} \frac{\alpha}{\gamma \tau_{j}^{dis*}+u}-\frac{\alpha \tau_{j}^{cen*}}{\zeta \tau_{j}^{cen*}+u}\right)^2\sigma_{W_j}^2}_{\text{(\ref{eq:variance-opt}.a)}} +\nonumber\\
~ & \underbrace{\left(\frac{\gamma-\zeta}{\zeta + u} \frac{\alpha}{\gamma \tau_{l}^{dis*}+u}+\frac{\alpha \tau_{l}^{dis*}}{\gamma \tau_{l}^{dis*}+u}-\frac{\alpha \tau_{l}^{cen*}}{\zeta \tau_{l}^{cen*}+u}\right)^2 \sigma_{W_l}^2}_{\text{(\ref{eq:variance-opt}.b)}}.
\end{align}
The term (\ref{eq:variance-opt}.a) equals to
\begin{align}
\left(\frac{\gamma+u}{\zeta + u} \frac{\alpha}{u}\right)^2\sigma_{W_j}^2, ~ & \mbox{when} \;  0 \le \sigma_{W_j}^2 < \bar{\tau}_1 \nonumber\\
\frac{\left(\alpha \gamma \sigma_{W_j}^2+u\sqrt{\frac{2mu}{IL}}\right)^2}{(\zeta+u)^2u^2\sigma_{W_j}^2}, ~ & \mbox{when} \; \sigma_{W_j}^2 \in [\bar{\tau}_1,\bar{\tau}_2) \nonumber\\
\frac{(\gamma+\frac{u}{\sqrt{L}})^2}{(\zeta+u)^2} \frac{2m}{Iu\sigma_{W_j}^2}, ~ & \mbox{when} \;  \sigma_{W_j}^2 \ge \bar{\tau}_2.\nonumber
\end{align}
Furthermore, the term (\ref{eq:variance-opt}.b) equals to
\begin{align}
\left(\frac{\gamma-\zeta}{\zeta + u} \frac{\alpha}{u}\right)^2\sigma_{W_l}^2 ,~ & \mbox{when} \; 	0 \le \sigma_{W_l}^2 < \bar{\tau}_1 \nonumber\\
\frac{\left(\alpha(\gamma-\zeta-u)\sigma_{W_l}^2+u\sqrt{\frac{2mu}{IL}}\right)^2}{(\zeta+u)^2u^2\sigma_{W_l}^2},~ & \mbox{when} \; \sigma_{W_l}^2 \in [\bar{\tau}_1,\bar{\tau}_2) \nonumber\\
\frac{[\gamma-\zeta+u(\sqrt{1/L}-1)]^2}{(\zeta+u)^2}\frac{2m}{Iu\sigma_{W_l}^2},~ & \mbox{when} \; \sigma_{W_l}^2 \ge \bar{\tau}_2 \nonumber
\end{align}
where $\bar{\tau}_1:=\sqrt{2mu/(IL)}/\alpha$, and $\bar{\tau}_2:=\sqrt{2mu/(I)}/\alpha$.

For any $\sigma_{W_j}$ in the first and second categories, the boundedness of $\sigma_{W_j}^2$ implies the boundedness of the corresponding term in \eqref{eq:variance-opt}. For any  $\sigma_{W_j}$ in the third category, the corresponding term in \eqref{eq:variance-opt} decreases with increasing $\sigma_{W_j}$ and is therefore also bounded. Thus, the variance $\mathbbm{D}(D_l^{cen*}-\hat{D}_l^{dis*})$ is bounded.

%\setcounter{equation}{0}  
%\renewcommand{\theequation}{F.\arabic{equation}}
%\subsection{Expected total profit}\label{apen-profit}
%\begin{align}
% ~ &   \mathbbm{E}[\sum_{l=1}^L \sum_{i \in \mathcal{I}_l} U_i(D_{i})-\lambda(D,W)D]-{\color{blue} \sum_{l=1}^L \sum_{i \in \mathcal{I}_l} h_l(\tau_{il})} \nonumber\\
% =~& \mathbbm{E}[tD-\frac{u}{2}\sum_{l=1}^L \sum_{i \in \mathcal{I}_l} D_i^2 -(\alpha D- \alpha \sum_{l=1}^L W_l + \beta_0) D]-\sum_{l=1}^L \frac{m}{\sigma_{W_l}^2}\frac{\tau_{l}}{1-\tau_{l}} \nonumber\\
% =~ & (t-\beta_0+\alpha \sum_{l=1}^L \overline{w}_l)IL b_1 -\frac{uIL}{2}(b_1^2+\sum_{l=1}^L (b_2^l)^2 (\sigma_{W_l}^2+\sigma_{\epsilon}^2)) \nonumber\\
% ~& -\alpha I^2 L^2 (b_1^2+\sum_{l=1}^L (b_2^l)^2 (\sigma_{W_l}^2+\sigma_{\epsilon}^2))-\alpha \sum_{l=1}^L (b_2^l)^2 \sigma_{W_l}^2 \nonumber\\
% ~& -\sum_{l=1}^L \sum_{i \in \mathcal{I}_l} h_l(\tau_{il})
%\end{align}	
	
\end{document}